# ON WEIGHTED *U*-STATISTICS FOR STATIONARY PROCESSES

BY TAILEN HSING AND WEI BIAO WU

*Texas A&M University and University of Chicago*

A weighted $U$-statistic based on a random sample $X_1, \ldots, X_n$ has the form $U_n = \sum_{1 \leq i,j \leq n} w_{i-j} K(X_i, X_j)$, where $K$ is a fixed symmetric measurable function and the $w_i$ are symmetric weights. A large class of statistics can be expressed as weighted $U$-statistics or variations thereof. This paper establishes the asymptotic normality of $U_n$ when the sample observations come from a nonlinear time series and linear processes.

**1. Introduction.** Consider the causal process

$$(1) \qquad X_i = F(\ldots, \varepsilon_{i-1}, \varepsilon_i),$$

where the $\varepsilon_j$ are i.i.d. random elements. Clearly (1) is very general and represents a huge class of processes. In particular, it contains the linear process $X_i = \sum_{j=0}^{\infty} a_j \varepsilon_{i-j}$, where $a_j$ are square summable and $\varepsilon_j$ has mean 0 and finite variance, and many nonlinear processes (cf. Section 3) including the threshold AR (TAR) models [Tong (1990)], AR with conditional heteroscedasticity (ARCH) models [Engle (1982)], random coefficient AR (RCA) models [Nicholls and Quinn (1982)], and exponential AR (EAR) models [Haggan and Ozaki (1981)]. The main goal of this paper is to consider the asymptotic behavior of the following statistic:

$$U_n = \sum_{1 \leq i,j \leq n} H_{i,j}(X_i, X_j) := \sum_{1 \leq i,j \leq n} w_{i-j} K(X_i, X_j),$$

where $K$ is a fixed symmetric measurable function and the $w_i$ are symmetric constants. We refer to $U_n$ as a weighted $U$-statistic. The class of statistics that can be written in this form or variations of this form is clearly huge. For example, if $H_{i_1,i_2}(x_1,x_2) = [G(x_1) + G(x_2)]/2$, $n^{-1} U_n$ is the partial sum of $G(X_1), \ldots, G(X_n)$; if $H_{i_1,i_2}(x_1,x_2) = x_1 x_2 I(|i_1 - i_2| = k)$, then $(n-k)^{-1} U_n$ 2 T. HSING AND W. B. WUis the sample covariance function of lag $k$ in $\{X_i\}$; if $H_{i_1,i_2} = I(i_1 \neq i_2)K$ and $K$, respectively, for some fixed function $K$, then $U_n$ is a (nonnormalized) $U$- and $V$-statistic, respectively.

The study of asymptotic properties of the weighted or even the usual $U$-statistics is in general not straightforward. Hoeffding's decomposition [Hoeffding (1961)] provides a powerful tool for understanding the large-sample properties of $U$-statistics based on i.i.d. or even weakly dependent observations. See Randles and Wolfe (1979), Serfling (1980) and Lee (1990). For the i.i.d. case, a small number of papers consider the asymptotic properties of weighted $U$-statistics; recent references include O'Neil and Redner (1993), Major (1994) and Rifi and Utzet (2000). For weakly dependent processes, the results for $U$-statistics are typically developed under mixing conditions; examples of these can be found in Yoshihara (1976), Denker and Keller (1983, 1986) and a series of recent papers by Borovkova, Burton and Dehling (1999, 2001, 2002). Laws of large numbers for $U$-statistics of stationary and ergodic sequences were considered by Aaronson, Burton, Dehling, Gilat, Hill and Weiss (1996) and Borovkova, Burton and Dehling (1999). For long-memory processes, $U$-statistics and quadratic forms were considered by Dehling and Taqqu (1989, 1991), Ho and Hsing (1996), Giraitis and Taqqu (1997) and Giraitis, Taqqu and Terrin (1998), among others.

Using martingale-based techniques, we prove some general results for $U_n$ for processes satisfying (1) in a variety of short- and long-memory situations. Approaches based on martingales are very effective in dealing with asymptotic issues of stationary processes. See Woodroofe (1992), Ho and Hsing (1996, 1997), Wu and Mielniczuk (2002) and Wu (2003) for some recent developments, where certain open problems are dealt with. Wu and Woodroofe (2004) investigate approximations to sums of stationary and ergodic sequences by martingales. Based on such approximations, they obtain necessary and sufficient conditions for such sums to be asymptotically normal from the martingale central limit theorem. No mixing conditions will be involved and the results obtained are often nearly optimal.

Specifically, in Section 2, we will state two general central limit theorems for a stationary process $Y_{i,j}$, where $Y_{i,j}$ is measurable with respect to the $\sigma$-field generated by $\varepsilon_k, k \leq i \vee j$, where $i \vee j = \max(i,j)$. An example of $Y_{i,j}$ is $Y_{i,j} = K(X_i, X_j)$, but the realm of possibilities goes beyond that. In addition to the dependence of the process $Y_{i,j}$, the $w_i$ introduce another level of dependence in $U_n$. The two cases of $\sum_{i=0}^{\infty} |w_i| < \infty$ and $\sum_{i=0}^{\infty} |w_i| = \infty$ correspond to short- and long-range dependence, respectively, thereby entailing norming sequences of different orders of magnitude. We will address both cases.

In Section 3, we apply the results to nonlinear time series that are geometric moment-contracting. These are "short-memory" processes, which include a large class of processes mentioned in the beginning of this section, and also processes that do not satisfy any strong mixing conditions. In Sections 4 and 5, respectively, our general results are applied to short- and long-memory linear processes. In the long-memory case, we let $Y_{i,j}$ be the remainder of an ANOVA decomposition of $K(X_i, X_j)$. The resulting decomposition of $U_n$ is similar in spirit to Hoeffding's decomposition, and the asymptotic distribution of $U_n$ can be determined by identifying the dominant term(s) of the decomposition. In Sections 3 and 4, we also compare some of our results with related results in Borovkova, Burton and Dehling (2001). The two sets of results have overlapping but somewhat different ranges of applicability; we explain the differences and, where they overlap, we point out situations where our results work more effectively.

Detailed proofs are included in Section 6.

**2. Notation and main results.** Let $\varepsilon_i, i \in \mathbb{Z}$, be i.i.d. random elements taking values in a general state space. Define the shift processes $\mathbf{Z}_i = (\ldots, \varepsilon_{i-1}, \varepsilon_i)$ and, for each $\ell \geq 1$, $\tilde{\mathbf{Z}}_i = \tilde{\mathbf{Z}}_{i,\ell} = (\varepsilon_{i-\ell+1}, \ldots, \varepsilon_i)$, where we often suppress $\ell$ in $\tilde{\mathbf{Z}}_{i,\ell}$ to simplify notation. Let $Y_{i,j}, i, j \in \mathbb{Z}$, be random variables with zero means and finite variances, such that $Y_{i,j} = Y_{j,i}$, $Y_{i,j} \in \sigma(\mathbf{Z}_{i \vee j})$ and $(Y_{i,j}, \mathbf{Z}_k)$ is a stationary process in the sense that the $(Y_{i+t,j+t}, \mathbf{Z}_{k+t})$ have the same finite-dimensional distributions as $(Y_{i,j}, \mathbf{Z}_k)$ for each $t \in \mathbb{Z}$; similarly let $\tilde{Y}_{i,j}, i, j \in \mathbb{Z}$, be random variables with zero means and finite variances, such that $\tilde{Y}_{i,j} = \tilde{Y}_{j,i}$, $\tilde{Y}_{i,j} \in \sigma(\tilde{\mathbf{Z}}_i, \tilde{\mathbf{Z}}_j)$ and $(\tilde{Y}_{i,j}, \tilde{\mathbf{Z}}_k)$ is a stationary process in the sense that the $(\tilde{Y}_{i+t,j+t}, \tilde{\mathbf{Z}}_{k+t})$ have the same finite-dimensional distributions as $(\tilde{Y}_{i,j}, \tilde{\mathbf{Z}}_k)$ for each $t \in \mathbb{Z}$. Define the projection operator

$$\mathcal{P}_t \xi = E(\xi | \mathbf{Z}_t) - E(\xi | \mathbf{Z}_{t-1}), \qquad t \in \mathbb{Z},$$

where $\xi$ is an integrable random variable. Let

$$L_{i,j} = w_{i-j} Y_{i,j} \quad \text{and} \quad \tilde{L}_{i,j} = w_{i-j} \tilde{Y}_{i,j}.$$

The two cases where the weights $w_i$ are summable and nonsummable have distinct flavors, and they will be considered separately in Sections 2.1 and 2.2.

2.1. *Summable weights.* In this section, we consider the asymptotic distribution of $\sum_{1 \leq i,j \leq n} w_{i-j} Y_{i,j}$, where the weights $w_i$ are absolutely summable. When $Y_{i,j} = K(X_i, X_j) - EK(X_i, X_j)$, obvious examples of this include partial sums for which $w_i = \delta_{i,0}$ and $k$-lag sample covariance function for which $w_i = \delta_{i,k}$. Let $\xrightarrow{d}$ denote convergence in distribution and let $\mathrm{N}(0, \sigma^2)$ be the normal distribution with mean zero and variance $\sigma^2$.

For any integers $i, j$, define

(2) $$\theta_{i,j} = \|\mathcal{P}_0 Y_{i,j}\|.$$

THEOREM 1. *Assume that*

$$\sum_{k=0}^{\infty}\sum_{i=0}^{\infty}|w_k|\theta_{i,i-k}<\infty. \tag{3}$$

*Then*

$$\frac{1}{\sqrt{n}}\sum_{1\leq i,j\leq n}L_{i,j}\xrightarrow{d}\mathrm{N}(0,\sigma^2) \tag{4}$$

*for some $\sigma^2<\infty$.*

REMARK 1. Since $\theta_{i,j}=\theta_{j,i}$ and $w_k=w_{-k}$, (3) is equivalent to the seemingly stronger statement

$$\sum_{k=-\infty}^{\infty}\sum_{i=0}^{\infty}|w_k|\theta_{i,i-k}<\infty \tag{5}$$

in view of

$$\sum_{k=-\infty}^{-1}\sum_{i=0}^{\infty}|w_k|\theta_{i,i-k}=\sum_{k=1}^{\infty}|w_k|\sum_{j=0}^{\infty}\theta_{j+k,j}$$

$$=\sum_{k=1}^{\infty}|w_k|\sum_{i=k}^{\infty}\theta_{i,i-k}\leq\sum_{k=0}^{\infty}\sum_{i=0}^{\infty}|w_k|\theta_{i,i-k}.$$

2.2. *Nonsummable weights.* The derivation of the main result, Theorem 3, in this section for nonsummable weights relies on Theorem 2, which asserts that $\sum_{1\leq i,j\leq n}L_{i,j}$ can be approximated by $\sum_{1\leq i,j\leq n}\tilde{L}_{i,j}$. To consider the asymptotic behavior of the latter, we apply the idea of the Hoeffding decomposition.

Let

$$\hat{\theta}_{i,j}=\sup_{\ell\geq 1}\|\mathcal{P}_0\tilde{Y}_{i,j}\|$$

and

$$\delta_\ell:=\sup_{j\in\mathbb{Z}}\|Y_{1,j}-\tilde{Y}_{1,j}\|. \tag{6}$$

Define

$$W_n(i)=\sum_{j=1}^{n}w_{i-j}\quad\text{and}\quad W_n=\left[\sum_{i=1}^{n}W_n^2(i)/n\right]^{1/2}.$$

THEOREM 2. *Assume that* $\lim_{\ell \to \infty} \delta_\ell = 0$, $\liminf_{n \to \infty} W_n / (\sum_{i=0}^{n} |w_i|) > 0$ *and*

$$\limsup_{\epsilon \to 0} \sum_{k \geq 0} \sum_{i=0}^{\infty} \min(\hat{\theta}_{i, i-k}, \epsilon) = 0. \tag{7}$$

*Then*

$$\lim_{\ell \to \infty} \limsup_{n \to \infty} \frac{1}{nW_n^2} \left\| \sum_{1 \leq i, j \leq n} (L_{i,j} - \tilde{L}_{i,j}) \right\|^2 = 0. \tag{8}$$

THEOREM 3. *Assume that* $\sum_{i=1}^{\infty} |w_i| = \infty$ *and* $\sum_{k=0}^{n}(n-k)w_k^2 = o(nW_n^2)$. *Then under the conditions of Theorem 2,*

$$\frac{1}{\sqrt{nW_n^2}} \sum_{1 \leq i, j \leq n} L_{i,j} \xrightarrow{d} \mathrm{N}(0, \sigma^2) \tag{9}$$

*for some* $\sigma^2 < \infty$.

REMARK 2. The assumptions on the $w_i$ in Theorems 2 and 3 are very minor and are satisfied for every situation of practical interest. For example, if $w_n \sim C/n^\beta$, $\beta < 1$, then those conditions hold. Note, however, that in Theorem 3, the second condition $\sum_{k=0}^{n}(n-k)w_k^2 = o(nW_n^2)$ cannot be derived from the first one $\sum_{i=1}^{\infty} |w_i| = \infty$. For example, let $w_n = 2^k$ whenever $n = 2^{2^k}$, $k \in \mathbb{N}$, and $w_n = 0$ otherwise. Then $\sum_{i=1}^{\infty} |w_i| = \infty$ and there exists a constant $c' > 0$ such that $\sum_{k=0}^{n}(n-k)w_k^2 \geq c'nW_n^2$ for all $n \geq 4$.

The conditions (3) and (7) are closely related through $\delta_\ell$. The following is useful in verifying the conditions in certain situations.

PROPOSITION 4. *The following hold*:

$$\sup_k \sum_{i=0}^{\infty} \theta_{i, i-k} \leq 2 \sum_{i=0}^{\infty} \delta_i \tag{10}$$

*and, for any* $\epsilon$,

$$\sup_k \sum_{i=0}^{\infty} \min(\hat{\theta}_{i, i-k}, \epsilon) \leq 4 \sum_{i=0}^{\infty} \min\left( \sup_{\ell \geq i} \delta_\ell, \epsilon \right). \tag{11}$$

PROOF. Let $j \geq i \geq \ell \geq 0$; then $\tilde{\mathbf{Z}}_i$ and $\tilde{\mathbf{Z}}_j$ are independent of $\mathbf{Z}_0$. Thus $\tilde{Y}_{i,j}$ is also independent of $\mathbf{Z}_0$ and $\mathcal{P}_0 \tilde{Y}_{i,j} = 0$. If $i \geq \ell, j \leq -1$, then $\tilde{\mathbf{Z}}_i$ is independent of $\mathbf{Z}_0$ and $\tilde{\mathbf{Z}}_j$ is $\mathbf{Z}_{-1}$ measurable. So $E[\tilde{Y}_{i,j}|\mathbf{Z}_0] = E[\tilde{Y}_{i,j}|\mathbf{Z}_{-1}]$ and, again, $\mathcal{P}_0 \tilde{Y}_{i,j} = 0$. Therefore,

$$
\begin{aligned}
\theta_{\ell,\ell-k} &= \|\mathcal{P}_0 Y_{\ell,\ell-k}\| = \|\mathcal{P}_0(Y_{\ell,\ell-k} - \tilde{Y}_{\ell,\ell-k})\| \\
&\leq \|Y_{\ell,\ell-k} - \tilde{Y}_{\ell,\ell-k}\| \leq \delta_\ell, \qquad k > \ell \geq 0, \\
\theta_{\ell+k,\ell} &= \|\mathcal{P}_0 Y_{\ell+k,\ell}\| = \|\mathcal{P}_0(Y_{\ell+k,\ell} - \tilde{Y}_{\ell+k,\ell})\| \\
&\leq \|Y_{\ell+k,\ell} - \tilde{Y}_{\ell+k,\ell}\| \leq \delta_\ell, \qquad k,\ell \geq 0,
\end{aligned}
\tag{12}
$$

by Cauchy's inequality. Hence

$$\sum_{i=0}^\infty \theta_{i,i-k} = \sum_{i=0}^{k-1} \theta_{i,i-k} + \sum_{i=0}^\infty \theta_{i+k,i} \leq 2 \sum_{i=0}^\infty \delta_i,$$

proving (10). To prove (11), similarly write

$$
\begin{aligned}
&\sum_{i=0}^\infty \min(\hat{\theta}_{i,i-k}, \epsilon) \\
&= \sum_{i=0}^{k-1} \min\left(\sup_{\ell \geq 0} \|\mathcal{P}_0 \tilde{Y}_{i,i-k}\|, \epsilon\right) + \sum_{i=0}^\infty \min\left(\sup_{\ell \geq 0} \|\mathcal{P}_0 \tilde{Y}_{i+k,i}\|, \epsilon\right) \\
&= \sum_{i=0}^{k-1} \min\left(\sup_{\ell \geq i} \|\mathcal{P}_0 \tilde{Y}_{i,i-k}\|, \epsilon\right) + \sum_{i=0}^\infty \min\left(\sup_{\ell \geq i} \|\mathcal{P}_0 \tilde{Y}_{i+k,i}\|, \epsilon\right).
\end{aligned}
\tag{13}
$$

Now, for $0 \leq i \leq k-1$, by the triangle inequality and (12),

$$\|\mathcal{P}_0 \tilde{Y}_{i,i-k}\| \leq \|\mathcal{P}_0 Y_{i,i-k}\| + \|\mathcal{P}_0(Y_{i,i-k} - \tilde{Y}_{i,i-k})\| \leq \delta_i + \delta_\ell,$$

where the same bound holds for $\|\mathcal{P}_0 \tilde{Y}_{i+k,i}\|$ if $i \geq 0$. Applying this and the inequality $\min(a+b,c) \leq \min(a,c) + \min(b,c)$ for $a,b,c \geq 0$, (11) follows readily from (13). $\square$

It follows from Proposition 4 that if both the $|w_i|$ and the $\delta_\ell$ are summable, then (3) holds; if $\sum_{i=0}^\infty \sup_{\ell \geq i} \delta_\ell < \infty$, then (7) holds.

**3. Nonlinear time series.** Let $\{\varepsilon'_j\}$ be an i.i.d. copy of $\{\varepsilon_j\}$. We say that $X_n = F(\mathbf{Z}_n)$ is *geometric moment-contracting* if there exist $\alpha > 0$, $C = C(\alpha) > 0$ and $0 < r(\alpha) < 1$ such that

$$
\begin{aligned}
E\{|F(\ldots, \varepsilon_{-1}, \varepsilon_0, \varepsilon_1, \ldots, \varepsilon_n) \\
- F(\ldots, \varepsilon'_{-1}, \varepsilon'_0, \varepsilon_1, \ldots, \varepsilon_n)|^\alpha\} \leq C r^n(\alpha), \qquad n \in \mathbb{N}.
\end{aligned}
\tag{14}
$$

Without loss of generality, let $\alpha < 1$ since otherwise we can employ the Hölder inequality. We may view $X'_n := F(\ldots, \varepsilon'_{-1}, \varepsilon'_0, \varepsilon_1, \ldots, \varepsilon_n)$ as a coupled version of $X_n$.

Condition (14) is very mild, and is satisfied by a wide class of nonlinear time series. Note that geometric moment contraction does not even require mixing (see Example 1). An important special class of (1) is the so-called *iterated random functions* such that (14) is satisfied. Let $X_n$ be defined recursively by

$$(15) \qquad X_n = F(X_{n-1}, \varepsilon_n),$$

where $F(\cdot, \cdot)$ is a bivariate measurable function with the Lipschitz constant

$$(16) \qquad L_\varepsilon = \sup_{x' \neq x} \frac{|F(x, \varepsilon) - F(x', \varepsilon)|}{|x - x'|} \leq \infty$$

satisfying

$$(17) \qquad E(\log L_\varepsilon) < 0.$$

Then the Markov chain (15) admits a unique stationary distribution if $E(L_\varepsilon^\alpha) < \infty$ and $E[|x_0 - F(x_0, \varepsilon)|^\alpha] < \infty$ for some $\alpha > 0$ and $x_0$ [Diaconis and Freedman (1999)]. The same set of conditions actually also imply the geometric-moment contraction (14) [cf. Lemma 3 in Wu and Woodroofe (2000)]. The condition (17) indicates that the iterated random function (15) contracts on average, which is satisfied for many popular nonlinear time series models such as TAR, RCA, ARCH and EAR under suitable conditions on model parameters.

Recall that $\mathbf{Z}_k = (\ldots, \varepsilon_{k-1}, \varepsilon_k)$ and $\tilde{\mathbf{Z}}_{k,\ell} = (\varepsilon_{k-\ell+1}, \ldots, \varepsilon_k)$. Let $\mathbf{Z}'_k = (\ldots, \varepsilon'_{k-1}, \varepsilon'_k)$.

LEMMA 5. *The geometric moment-contraction condition* (14) *holds if and only if there exist* $F_1, F_2, \ldots$, *with each* $F_\ell$ *being an $\ell$-variate measurable function, such that, for some* $C < \infty$,

$$(18) \qquad E\{|F(\mathbf{Z}_k) - F_\ell(\tilde{\mathbf{Z}}_{k,\ell})|^\alpha\} \leq Cr^\ell(\alpha), \qquad \ell \in \mathbb{N}.$$

PROOF. The "$\Rightarrow$" direction. Assume (14). Then for each $\ell$, there exists a realization $\mathbf{Z}'_0 = \mathbf{z}_0$ such that $E(|X_\ell - X'_\ell|^\alpha | \mathbf{Z}'_0 = \mathbf{z}_0) \leq Cr^\ell(\alpha)$. So (18) holds by defining $F_\ell(\cdot) = F(\mathbf{z}_0, \cdot)$, which is clearly measurable. The "$\Leftarrow$" direction follows easily from

$$\begin{aligned}
E(|X_\ell - X'_\ell|^\alpha) &= E[|F(\mathbf{Z}_k) - F(\mathbf{Z}'_{k-\ell}, \tilde{\mathbf{Z}}_{k,\ell})|^\alpha] \\
&\leq E[|F(\mathbf{Z}_k) - F_\ell(\tilde{\mathbf{Z}}_{k,\ell})| + |F(\mathbf{Z}'_{k-\ell}, \tilde{\mathbf{Z}}_{k,\ell}) - F_\ell(\tilde{\mathbf{Z}}_{k,\ell})|]^\alpha \\
&\leq E[|F(\mathbf{Z}_k) - F_\ell(\tilde{\mathbf{Z}}_{k,\ell})|^\alpha] + E[|F(\mathbf{Z}'_{k-\ell}, \tilde{\mathbf{Z}}_{k,\ell}) - F_\ell(\tilde{\mathbf{Z}}_{k,\ell})|^\alpha] \\
&= 2E[|F(\mathbf{Z}_k) - F_\ell(\tilde{\mathbf{Z}}_{k,\ell})|^\alpha],
\end{aligned}$$

where we have applied the inequality $|a + b|^\alpha \leq |a|^\alpha + |b|^\alpha$ for $0 < \alpha \leq 1$. □

In Lemma 5, we can often choose $\mathbf{z}_0$ arbitrarily in defining $F_\ell$. This can be illustrated by the correlation integral example in Theorem 7.

We remark that conditions similar to (14) and (18) have appeared in the literature. Denker and Keller (1986) assumed that $F$ is Lipschitz-continuous in the sense that there exists a $\rho \in (0,1)$ for which

$$(19) \qquad |F(\ldots, z_{n-1}, z_n) - F(\ldots, z'_{n-1}, z'_n)| \leq \text{const.} \cdot \rho^n$$

if $z_1 = z'_1, \ldots, z_n = z'_n$. For the two-sided extension, see Definition 1.3 in Borovkova, Burton and Dehling (2001). Comparing with our condition (14), (19) is stronger and it does not allow models like $X_n = \rho X_{n-1} + \varepsilon_n$, where $|\rho| < 1$ and the random variables $\varepsilon_n$ are i.i.d. with unbounded support. Borovkova, Burton and Dehling (2001) proposed a weaker version of (19), termed *r-approximation condition*, which requires

$$(20) \qquad d_l(r) := E|X_0 - E(X_0|\varepsilon_{-l}, \ldots, \varepsilon_l)|^r \to 0 \qquad \text{as } \ell \to \infty,$$

for some $r \geq 1$. To make this weaker version operational, one needs to implicitly assume $E|X_0| < \infty$, which excludes the case that $\varepsilon_n$ does not have a mean. Our formulation has the advantage that heavy-tailed distributions are allowed.

As before, write $X_i = F(\mathbf{Z}_i)$, and for a fixed choice of $\ell$-variate function $F_\ell$ from Lemma 5, define $\tilde{X}_i = \tilde{X}_{i,\ell} = F_\ell(\tilde{\mathbf{Z}}_{i,\ell})$; let $Y_{i_1,i_2} = K(X_{i_1}, X_{i_2}) - E[K(X_{i_1}, X_{i_2})]$ and $\tilde{Y}_{i_1,i_2} = K(\tilde{X}_{i_1}, \tilde{X}_{i_2}) - EK(\tilde{X}_{i_1}, \tilde{X}_{i_2})$ and recall that $L_{i_1,i_2} = w_{i_1-i_2} Y_{i_1,i_2}$ and $\tilde{L}_{i_1,i_2} = w_{i_1-i_2} \tilde{Y}_{i_1,i_2}$. Then Theorems 1 and 3 imply (i) and (ii) of the following result, respectively, in view of Proposition 4.

THEOREM 6. *Suppose that for each $\ell \geq 1$, there exists an $\ell$-variate function $F_\ell$ such that (18) holds and $\sum_{i=0}^\infty \sup_{\ell \geq i} \delta_\ell < \infty$.*

(i) *If $\sum |w_i| < \infty$, then*

$$n^{-1/2} \sum_{i,j=1}^n w_{i-j}[K(X_i, X_j) - EK(X_i, X_j)] \xrightarrow{d} N(0, \sigma^2)$$

*for some $\sigma^2 < \infty$.*

(ii) *Let $\sum_{i=1}^\infty |w_i| = \infty$. If we also have $\liminf_{n \to \infty} W_n/(\sum_{i=0}^n |w_i|) > 0$ and $\sum_{k=0}^n (n-k) w_k^2 = o(nW_n^2)$, then*

$$(nW_n^2)^{-1/2} \sum_{i,j=1}^n w_{i-j}[K(X_i, X_j) - EK(X_i, X_j)] \xrightarrow{d} N(0, \sigma^2)$$

*for some $\sigma^2 < \infty$.*

The inequality (18) implies that the distance between $X_i$ and $\tilde{X}_{i,\ell}$ decays exponentially fast to 0 in $\ell$. Thus under certain continuity conditions on $K$, $\delta_\ell$ is expected to vanish sufficiently quickly. For an application of Theorem 6, consider the correlation integral

$$N_b = \sum_{i_1,i_2=1}^n \mathbf{1}_{|X_{i_1}-X_{i_2}|<b},$$

which measures the number of pairs $(X_i, X_j)$ such that their distance is less than $b > 0$. Correlation integral is of critical importance in the study of dynamical systems; see Wolff (1990), Serinko (1994) and Denker and Keller (1986) for further references.

THEOREM 7. *Suppose that $X_n$ defined in* (1) *satisfies* (14), *and for some $\kappa > 1$,*

(21) $$\sup_{j\neq 0,\, x\in \mathbf{R}} P(x < X_0 - X_j \le x + \tau) \le C \log^{-2\kappa} \tau^{-1}$$

*for all $0 < \tau < 1/2$. Then $[N_b - E(N_b)]/n^{3/2} \xrightarrow{d} N(0, \sigma^2)$ for some $\sigma^2 < \infty$.*

PROOF. Let $K(x,y) = \mathbf{1}_{|x-y|<b}$ and $w_i \equiv 1$. By Lemma 5, for each $\ell$ there exists an $\ell$-variate measurable function $F_\ell(\cdot)$ such that (18) holds. Next we shall verify that (18) together with (21) implies that $\delta_\ell = O(\ell^{-\kappa})$, which is summable and thus completes the proof in view of (ii) of Theorem 6. To this end, let $u = r(\alpha)^{1/(2\alpha)} < 1$. Then by (18) and the Markov inequality, $\|\mathbf{1}_{|X_0-\tilde{X}_0|\ge u^\ell}\|^2 \le u^{-\alpha\ell} E|X_0 - \tilde{X}_0|^\alpha \le C r(\alpha)^{\ell/2}$, where, as usual, $\tilde{X}_i = \tilde{X}_{i,\ell}$. For any $0 < u < 1$, observe that

$$\|[K(X_0, X_i) - K(\tilde{X}_0, \tilde{X}_i)]\mathbf{1}_{\max(|X_0-\tilde{X}_0|,|X_i-\tilde{X}_i|)<u^\ell}\|$$
$$\le P^{1/2}(b - 2u^\ell \le |X_0 - X_i| \le b) + P^{1/2}(b \le |X_0 - X_i| \le b + 2u^\ell)$$

which, by (21), is bounded by $2C^{1/2} \log^{-\kappa}(2u^\ell)^{-1} = O(\ell^{-\kappa})$. Since $|K| \le 1$,

$$\|[K(X_0, X_i) - K(\tilde{X}_0, \tilde{X}_i)]\|$$
$$\le \|[K(X_0, X_i) - K(\tilde{X}_0, \tilde{X}_i)]\mathbf{1}_{\max(|X_0-\tilde{X}_0|,|X_i-\tilde{X}_i|)<u^\ell}\|$$
$$+ \|[K(X_0, X_i) - K(\tilde{X}_0, \tilde{X}_i)]\mathbf{1}_{\max(|X_0-\tilde{X}_0|,|X_i-\tilde{X}_i|)\ge u^\ell}\|$$
$$= O(\ell^{-\kappa}) + O[r(\alpha)^{\ell/4}],$$

proving that $\delta_\ell = O(\ell^{-\kappa})$. □

EXAMPLE 1. Let $X_n = (X_{n-1} + \varepsilon_n)/2$, where $\varepsilon_n$ are i.i.d. Bernoulli random variables with success probability $1/2$. Then $X_n$ admits Uniform$(0,1)$ as a stationary distribution. This process is not strong mixing. Now we show that (21) is satisfied. Assume $j \geq 1$. Let $U = \sum_{i=1}^{j} \varepsilon_i/2^{j-i}$. Then $U$ is uniformly distributed over $\{0, 1/2^j, \ldots, (2^j-1)/2^j\}$ and $X_j = X_0/2^j + U$. Hence (21) holds in view of

$$\begin{aligned}
P(x < X_0 - X_j \leq x + \tau) \\
= EP[(x+U)/(1-2^{-j}) < X_0 \leq (x+\tau+U)/(1-2^{-j})|U] \\
\leq \tau/(1-2^{-j}) \leq 2\tau.
\end{aligned}$$

The process $X_n$ is related to the doubling map $Tx := 2x \bmod 1$ in the following way. Let $Y_0$ be a Uniform$(0,1)$ random variable and define recursively $Y_i = 2Y_{i-1} \bmod 1$ for $i \geq 1$. Then $(X_1, \ldots, X_n)$ has the same distribution as $(Y_n, \ldots, Y_1)$ and hence $N_b$ and $M_b = \sum_{i,j=1}^{n} \mathbf{1}_{|Y_i - Y_j| < b}$ are identically distributed. The limiting distribution of the empirical $U$-process $\{M_b, 0 \leq b \leq 1\}$ was discussed in Borovkova, Burton and Dehling (2001); see Section 6 therein.

Our Theorems 6 and 7 are closely related to certain results by Borovkova, Burton and Dehling (2001), which considered nonweighted $U$-statistics for two-sided processes $X_n = F((\varepsilon_{n+k})_{k \in \mathbb{Z}})$, where $(\varepsilon_n)_{n \in \mathbb{Z}}$ is stationary and absolutely regular (or weak Bernoulli). To make a specific comparison, we state here their Theorem 7, a central limit theorem. Let

$$\beta_k = 2\sup_n\{\sup\{P(A|\varepsilon_1, \ldots, \varepsilon_n) - P(A) : A \text{ is } \sigma(\varepsilon_{n+k}, \varepsilon_{n+k+1}, \ldots)\text{-measurable}\}\}$$

be the mixing coefficients, let $\alpha_k = [2\sum_{i=k}^{\infty} d_k(1)]^{1/2}$ with $d_k(1)$ defined by (20), and let $K(\cdot, \cdot)$ be a bounded, symmetric function such that

$$\sup_{1 \leq k \leq \infty} E\{|K(X_0, X_k) - K(X', X_k)|\mathbf{1}_{|X_0 - X'| \leq \tau}\} \leq \phi(\tau)$$

with $\lim_{\tau \to 0} \phi(\tau) = 0$, where $X'$ is identically distributed as $X_0$ and $X_\infty$ is interpreted as an independent copy of $X_0$. See Definitions 1.2, 1.4 and 2.12 in Borovkova, Burton and Dehling (2001). Then the asymptotic normality of $U_n$ holds provided

$$\sum_{k=1}^{\infty} k^2(\beta_k + \alpha_k + \phi(\alpha_k)) < \infty. \tag{22}$$

This result has a number of similarities with our Theorem 6. However, the two results do not imply one another. Theorem 6 assumes one-sided processes with i.i.d. innovations while their result allows the innovations to be two-sided and weakly dependent; on the other hand, Theorem 6 allows unbounded $K$, general weights $w_i$ and process $X_k$ for which the mean is infinite, whereas their result requires $K$ to be bounded, $w_i = 1$ and $E|X_0| <$ $\infty$. Let us make a more specific comparison in the context of Theorem 7 of the present paper where both results are applicable. Applying the central limit theorem in Borovkova, Burton and Dehling (2001) to $N_b$ for one-sided processes with i.i.d. innovations satisfying (18) with $\alpha = 1$ and (21), we have $\beta_n = 0$ and $\alpha_n \leq C\rho^n$ for all $n \geq 1$, where $\rho \in (0,1)$. By Example 2.2 in Borovkova, Burton and Dehling (2001),

$$\phi(\tau) = \sup_{1 \leq k \leq \infty} P(b - \tau \leq |X_0 - X_k| \leq b + \tau) \leq 2C \log^{-2\kappa} \tau^{-1}.$$

Thus condition (22) is reduced to $\sum_{k=1}^{\infty} k^2 \log^{-2\kappa}(1/\rho^k) < \infty$, namely $\kappa > 3/2$, which is stronger than the condition $\kappa > 1$ imposed in Theorem 7.

**4. Short-memory linear processes.** Let $(a_n)_{n \geq 0}$ be square summable, let $(\varepsilon_n)_{n \geq 0}$ be i.i.d. random variables with mean 0 and finite variance and let

$$X_n = \sum_{i=0}^{\infty} a_i \varepsilon_{n-i}. \tag{23}$$

If $\sum_{i=0}^{\infty} |a_i| < \infty$, then the covariance function $\Gamma(n) = E(X_0 X_n)$ is summable and we say that $X_n$ is short-memory. In this section, let $Y_{i,j} = K(X_i, X_j) - EK(X_i, X_j)$. For short-memory processes, we shall utilize the linearity structure and provide conditions on $K(\cdot, \cdot)$ such that (3) and (7) hold by computing the quantities $\theta_{i,j}$ in (2) and $\hat{\theta}_{i,j}$ in (7). In Section 5, we shall discuss the case when $X_n$ is long-memory, which has a very different flavor. Note that $\sum_{\ell=0}^{\infty} \delta_\ell < \infty$ is not guaranteed by $\sum_{i=0}^{\infty} |a_i| < \infty$, and we need more refined computations, which is feasible by the linearity of $X_n$.

Let $\tilde{a}_i = a_i I(i < \ell)$, $\tilde{X}_n = \sum_{i=0}^{\infty} \tilde{a}_i \varepsilon_{n-i}$ and $\tilde{Y}_{i,j} = K(\tilde{X}_i, \tilde{X}_j) - EK(\tilde{X}_i, \tilde{X}_j)$. Also define $X'_{i,j_1,j_2} = X_{i,j_1,j_2} - a_i \varepsilon_0 + a_i \varepsilon'_0$, where the truncated process

$$X_{i,j_1,j_2} = \begin{cases} \sum_{j_1 \leq j \leq j_2} a_{i-j} \varepsilon_j, & -\infty \leq j_1 \leq j_2 \leq \infty, \\ 0, & -\infty \leq j_2 + 1 \leq j_1 \leq \infty. \end{cases}$$

Define the convolutions

$$\begin{aligned} K_{i_1,i_2,j}(x_1, x_2) &= EK(x_1 + X_{i_1,j+1,\infty}, x_2 + X_{i_2,j+1,\infty}), \\ K_{i_1,j}(x_1, x_2) &= EK(x_1 + X_{i_1,j+1,\infty}, x_2), \qquad x_1, x_2 \in \mathbf{R}. \end{aligned} \tag{24}$$

Let $\tilde{K}_{i_1,i_2,j}$, $\tilde{K}_{i_1,j}$, $\tilde{X}_{i,j_1,j_2}$ and $\tilde{X}'_{i,j_1,j_2}$ be defined similarly to $K_{i_1,i_2,j}$, $K_{i_1,j}$, $X_{i,j_1,j_2}$ and $X'_{i,j_1,j_2}$ with $\mathbf{Z}_i$ replaced by $\tilde{\mathbf{Z}}_i$.

PROPOSITION 8. *Assume that* $\sup_{i,j} \|K(X_i, X_j)\| < \infty$ *and* $\sup_{i,j,\ell} \|K(\tilde{X}_i, \tilde{X}_j)\| < \infty$. *Further, assume that there exist* $n_0 \in \mathbb{N}$ *and* $C < \infty$ *such that,*

*for all* $i_1, i_2, \ell \geq n_0$,

$$
\begin{aligned}
(25) \quad &\|\tilde{K}_{i_1,i_2,0}(\tilde{X}'_{i_1,-\infty,0}, \tilde{X}'_{i_2,-\infty,0}) - \tilde{K}_{i_1,i_2,0}(\tilde{X}_{i_1,-\infty,0}, \tilde{X}_{i_2,-\infty,0})\| \\
&\leq C(|a_{i_1}| + |a_{i_2}|)
\end{aligned}
$$

*and*

$$
(26) \quad \sup_{k \leq -1} \|\tilde{K}_{i_1,0}(\tilde{X}'_{i_1,-\infty,0}, X_k) - \tilde{K}_{i_1,0}(\tilde{X}_{i_1,-\infty,0}, X_k)\| \leq C|a_{i_1}|.
$$

*Then the following hold.*

(i)

$$
\hat{\theta}_{i,j} \leq \begin{cases} C(|a_i| + |a_j|), & i, j \geq n_0, \\ C|a_i|, & i \geq n_0, \ j < 0. \end{cases}
$$

(ii) *There exists some constant* $C' < \infty$ *such that*

$$
\sup_{k \geq 0} \sum_{i=0}^{\infty} \theta_{i,i-k} \leq C' \sum_{i=0}^{\infty} \sup_{j \geq i} |a_i|,
$$

*and for any* $\epsilon > 0$,

$$
\sup_{k \geq 0} \sum_{i=0}^{\infty} \min(\hat{\theta}_{i,i-k}, \epsilon) \leq C' \left[ \epsilon + \sum_{i=0}^{\infty} \min\left(\sup_{j \geq i} |a_i|, \epsilon\right) \right].
$$

PROOF. Fix $i_1, i_2 \geq n_0$. First we remark that if $\ell < n_0$, then the left-hand sides of (25) and (26) are both equal to 0 so that the inequalities trivially hold. Writing

$$
E[\tilde{K}_{i_1,i_2,0}(\tilde{X}'_{i_1,-\infty,0}, \tilde{X}'_{i_2,-\infty,0})|\mathbf{Z}_0] = \tilde{K}_{i_1,i_2,-1}(\tilde{X}_{i_1,-\infty,-1}, \tilde{X}_{i_2,-\infty,-1}),
$$

we have by Cauchy's inequality that

$$
\begin{aligned}
&\|\mathcal{P}_0 K(\tilde{X}_{i_1}, \tilde{X}_{i_2})\| \\
&= \|\tilde{K}_{i_1,i_2,0}(\tilde{X}_{i_1,-\infty,0}, \tilde{X}_{i_2,-\infty,0}) - \tilde{K}_{i_1,i_2,-1}(\tilde{X}_{i_1,-\infty,-1}, \tilde{X}_{i_2,-\infty,-1})\| \\
&= \|E[\tilde{K}_{i_1,i_2,0}(\tilde{X}_{i_1,-\infty,0}, \tilde{X}_{i_2,-\infty,0}) - \tilde{K}_{i_1,i_2,0}(\tilde{X}'_{i_1,-\infty,0}, \tilde{X}'_{i_2,-\infty,0})|\mathbf{Z}_0]\| \\
&\leq \|\tilde{K}_{i_1,i_2,0}(\tilde{X}'_{i_1,-\infty,0}, \tilde{X}'_{i_2,-\infty,0}) - \tilde{K}_{i_1,i_2,0}(\tilde{X}_{i_1,-\infty,0}, \tilde{X}_{i_2,-\infty,0})\|.
\end{aligned}
$$

Similarly, $\|\mathcal{P}_0 K(\tilde{X}_{i_1}, \tilde{X}_k)\| \leq \|\tilde{K}_{i_1,0}(\tilde{X}'_{i_1,-\infty,0}, X_k) - \tilde{K}_{i_1,0}(\tilde{X}_{i_1,-\infty,0}, X_k)\|$, which completes the proof of (i) in view of (25) and (26).

The first inequality in (ii) follows simply from

$$
\sum_{i=0}^{\infty} \theta_{i,i-k} = \sum_{i=0}^{k-1} \theta_{i,i-k} + \sum_{i=0}^{\infty} \theta_{i+k,i} \leq C' + 2\sum_{i=n_0}^{\infty} \sup_{j \geq i} |a_i|,
$$

whereas the second inequality there can be derived similarly. □

REMARK 3. Note that $\tilde{X}'_{i_1,-\infty,0} - \tilde{X}_{i_1,-\infty,-1} = a_{i_1}(\varepsilon'_0 - \varepsilon_0)$, conditions (25) and (26) can be interpreted as the "Lipschitz continuity" of $\tilde{K}_{i_1,i_2,0}$ and $\tilde{K}_{i_1,0}$. Discontinuous functions $K$ are allowed since these are convolutions of $K$ and the distribution functions of $(\tilde{X}_{i_1,1,\infty}, \tilde{X}_{i_2,1,\infty})$ and $(\tilde{X}_{i_1,j+1,\infty}, 0)$, respectively. For example, if $K$ is a bounded function and $f_\varepsilon$, the density function of $\varepsilon_1$, satisfies $\int |f'_\varepsilon(t)| \, dt < \infty$, then it is easily seen that (25) and (26) hold. Observe that degree of smoothness of the distributions of the above random vectors increases with $i_1, i_2$. Thus, by only requiring (25) and (26) to hold for large $i_1, i_2$, an additional dimension of flexibility is in place.

Proposition 8 together with Theorems 1 and 3 immediately yield

THEOREM 9. *Assume that $\sum_{i=1}^{\infty} \sup_{j \geq i} |a_j| < \infty$. Also assume that $\sup_{i,j} \|K(X_i, X_j)\| < \infty$ and $\sup_{i,j,\ell} \|K(\tilde{X}_i, \tilde{X}_j)\| < \infty$ and that the regularity conditions (25) and (26) hold.*

(i) *If $\sum_{i=0}^{\infty} |w_i| < \infty$, then*
$$n^{-1/2} \sum_{i,j=1}^{n} w_{i-j}[K(X_i, X_j) - EK(X_i, X_j)] \xrightarrow{d} \mathrm{N}(0, \sigma^2)$$
*for some $\sigma^2 < \infty$.*

(ii) *Suppose that $\sum_{i=1}^{\infty} |w_i| = \infty$ with $\liminf_{n \to \infty} \frac{W_n}{\sum_{i=0}^{n} |w_i|} > 0$ and $\sum_{k=0}^{n} (n-k)w_k^2 = o(nW_n^2)$. Assume also that $\lim_{\ell \to \infty} \delta_\ell = 0$. Then*
$$(nW_n^2)^{-1/2} \sum_{i,j=1}^{n} w_{i-j}[K(X_i, X_j) - EK(X_i, X_j)] \xrightarrow{d} \mathrm{N}(0, \sigma^2)$$
*for some $\sigma^2 < \infty$.*

REMARK 4. Consider the special case in which $a_n = n^{-\beta}$ for $n \geq 1$ and $\beta > 1$ and $\varepsilon_i$ are i.i.d. standard normal random variables. Then by (20), $d_n(1) \sim c_1 n^{1/2-\beta}$ for some $c_1 > 0$. Here $\gamma_n \sim \beta_n$ is meant as $\lim_{n \to \infty} \gamma_n / \beta_n = 1$. So $\alpha_n \sim c_2 n^{3/4-\beta/2}$ and the condition (22) of Borovkova, Burton and Dehling (2001) necessarily requires $\sum_{n=1}^{\infty} n^2 \alpha_n < \infty$, or $\beta > 15/2$. In comparison, our Theorem 9 only imposes $\beta > 1$.

**5. Long-memory linear processes.** In (23), let $a_j = j^{-\beta} L(j) I(j \geq 1)$ for some $\beta \in (1/2, 1)$ and slowly varying function $L$. Thus, $a_j$ is regularly varying at $\infty$ with index $-\beta$. This represents a rich class of processes. In particular, it contains the important time series model fractional autoregressive integrated moving average (FARIMA) process. See Granger and Joyeux (1980). Note that $\{X_i\}$ is long-range dependent in the sense that the covariances are not summable [cf. Beran (1994)].

Let $\mathcal{K} = \{0\} \cup \{\mathbf{k} = (k_1, \ldots, k_r) : r = 1, 2, \ldots, k_i \in \{1, 2\}\}$ and let $|\mathbf{k}|$ be the length of $\mathbf{k}$ ($|0| = 0$). We assume throughout the section that integer $\rho \geq 1$ satisfies

$$\tag{27} \sum_{n=1}^{\infty} n^{-\beta(\rho+1)+\rho/2} |L(n)|^{\rho+1} < \infty.$$

Condition (27) allows simultaneous consideration of two cases: (i) $(\rho+1)(2\beta-1) > 1$ and (ii) $(\rho+1)(2\beta-1) = 1$ and $\sum_{n=1}^{\infty} |L^{\rho+1}(n)|/n < \infty$. Case (i) has been widely studied [see Ho and Hsing (1997)], while the boundary case $(\rho+1)(2\beta-1) = 1$ has been overlooked in the literature. Our approach allows us to investigate the boundary case for which the limiting behavior depends on the growth of the slowly varying function $L$.

Denote by $\mathbf{C}^{\rho}(\mathbf{R}^2)$ the class of all functions $g$ such that the partial derivatives $D_\mathbf{k} g = \partial^r g / \partial x_{k_1} \cdots \partial x_{k_r}$ exist for all $\mathbf{k} = (k_1, \ldots, k_r) \in \mathcal{K}$ for which $|\mathbf{k}| \leq \rho$. For each $i_1, i_2$, let

$$Y_{i_1,i_2} = K(X_{i_1}, X_{i_2}) - \sum_{r=0}^{\rho} \sum_{|\mathbf{l}|=r} D_\mathbf{l} K_{i_1,i_2,-\infty}(0,0) \sum_{j_1 > \cdots > j_r} \prod_{s=1}^{r} a_{i_{l_s}-j_s} \varepsilon_{j_s}$$

and $L_{i_1,i_2} = w_{i_1-i_2} Y_{i_1,i_2}$. Let $\tilde{Y}_{i,j}, \tilde{L}_{i,j}$ be defined as $Y_{i,j}, L_{i,j}$ with $\tilde{a}_i = a_i I(i < \ell)$ replacing $a_i$. Let $K_{i_1,i_2,j}$, $K_{i_1,j}$, $X_{i,j_1,j_2}$, $\tilde{K}_{i_1,i_2,j}$, $\tilde{K}_{i_1,j}$ and $\tilde{X}_{i,j_1,j_2}$ be defined as in Section 4.

Write

$$\tag{28} U_n = \sum_{i_1,i_2=1}^{n} L_{i_1,i_2} + \sum_{r=0}^{\rho} Z_{n,r},$$

where

$$Z_{n,r} := \sum_{|\mathbf{l}|=r} \sum_{i_1,i_2=1}^{n} w_{i_1-i_2} D_\mathbf{l} K_{i_1,i_2,-\infty}(0,0) \sum_{j_1 > \cdots > j_r} \prod_{s=1}^{r} a_{i_{l_s}-j_s} \varepsilon_{j_s}.$$

Observe that $Z_{n,r}, 1 \leq r \leq \rho$, are well-structured, and can be shown to follow non-central limit theorems under mild regularity conditions on the $D_\mathbf{l} H_{i_1,i_2,-\infty}(0,0)$. Our main results, Theorems 10 and 11, show that the normalized $\sum_{i_1,i_2=1}^{n} L_{i_1,i_2}$ follows a central limit theorem under mild conditions. These two pieces of information will then combine to give a comprehensive picture of the asymptotic behavior of $U_n$. We refer to $\sum_{i_1,i_2=1}^{n} L_{i_1,i_2}$ and $\sum_{r=1}^{\rho} Z_{n,r}$, respectively, as the short- and long-memory components of $U_n$.

We now state the technical conditions for our main results. In the following, let

$$A_i(k) = \sum_{j=i}^{\infty} a_j^k, \qquad k=2,4, \ i \geq 0.$$

(K1) There exists $n_0 \in \mathbb{N}$ such that $K_{i_1,i_2,0}(\cdot,\cdot) \in \mathbf{C}^\rho(\mathbf{R}^2)$ when $i_1,i_2 \geq n_0$, and for all $\mathbf{k} \in \mathcal{K}$ with $|\mathbf{k}| \leq \rho$,

(29)
$$E[D_{\mathbf{k}} K_{i_1,i_2,0}(X_{i_1,-\infty,0}, X_{i_2,-\infty,0}) | \mathbf{Z}_{-1}]$$
$$= D_{\mathbf{k}} K_{i_1+1,i_2+1,0}(X_{i_1,-\infty,-1}, X_{i_2,-\infty,-1}).$$

(K2) For $i_1, i_2 \geq n_0$, there exists $C < \infty$ such that, for all $\mathbf{k} \in \mathcal{K}$ with $|\mathbf{k}| < \rho$,

(30)
$$\| D_{\mathbf{k}} K_{i_1,i_2,0}(X_{i_1,-\infty,0}, X_{i_2,-\infty,0}) - D_{\mathbf{k}} K_{i_1,i_2,-1}(X_{i_1,-\infty,-1}, X_{i_2,-\infty,-1})$$
$$- \langle \nabla D_{\mathbf{k}} K_{i_1,i_2,-1}(X_{i_1,-\infty,-1}, X_{i_2,-\infty,-1}), (a_{i_1}\varepsilon_0, a_{i_2}\varepsilon_0) \rangle \|$$
$$\leq C(a_{i_1}^2 + a_{i_2}^2),$$

and, for $|\mathbf{k}| = \rho$,

(31)
$$\| D_{\mathbf{k}} K_{i_1,i_2,0}(X_{i_1,-\infty,0}, X_{i_2,-\infty,0}) - D_{\mathbf{k}} K_{i_1,i_2,-1}(0,0) \|^2$$
$$\leq C[A_{i_1}(2) + A_{i_2}(2)].$$

(K3) For $i_1 \geq n_0$, there exists $C < \infty$ such that

(32)
$$\sup_{k \leq -1} \| \mathcal{P}_0 K(X_{i_1}, X_k) \|$$
$$= \sup_{k \leq -1} \| K_{i_1,1}(X_{i_1,-\infty,0}, X_k) - K_{i_1,0}(X_{i_1,-\infty,-1}, X_k) \| \leq C |a_{i_1}|.$$

Similarly, we define

(K4) There exists $n_0 \in \mathbb{N}$ such that, for each $y$, $K_{i,0}(\cdot, y) \in \mathbf{C}^\rho(\mathbf{R})$ when $i_1 \geq n_0, i_2 \leq -1$, and, for all $k \leq \rho$,

(33) $$E[K_{i_1,0}^{(k,0)}(X_{i_1,-\infty,0}, X_{i_2}) | \mathbf{Z}_{-1}] = K_{i_1+1,0}^{(k,0)}(X_{i_1,-\infty,-1}, X_{i_2}).$$

(K5) For $i_1 \geq n_0, i_2 \leq -1$, there exists $C < \infty$ such that, for all $k < \rho$,

(34)
$$\| K_{i_1,0}^{(k,0)}(X_{i_1,-\infty,0}, X_{i_2}) - K_{i_1,i_2,-1}^{(k)}(X_{i_1,-\infty,-1}, X_{i_2})$$
$$- K_{i_1,-1}^{(k+1,0)}(X_{i_1,-\infty,-1}, X_{i_2}), \ a_{i_1}\varepsilon_0 \| \leq C a_{i_1}^2,$$

and, for all $k \leq \rho$,

(35) $$\| K_{i_1,0}^{(k,0)}(X_{i_1,-\infty,0}, X_{i_2}) - K_{i_1,-1}^{(k,0)}(0, X_{i_2}) \|^2 \leq C A_{i_1}(2).$$

(K6)

(36) $$\lim_{\ell \to \infty} \sup_{j \geq 1} \|K(X_1, X_j) - K(\tilde{X}_1, \tilde{X}_j)\| = 0.$$

REMARK 5. The conditions (K1) and (K4) state that we can interchange the order of differentiation and integration. The other conditions are smoothness conditions on $K_{i_1,i_2,0}$ and $K_{i_1,0}$ for large $i_1, i_2$. For the latter, Remarks 3 and 4 are still relevant. These conditions are left in the form in which they are directly applied in the proofs. Finding sufficient conditions that are easy to work with in specific contexts should be straightforward. See Ho and Hsing (1997), Koul and Surgailis (1997) and Giraitis and Surgailis (1999).

In the following, we consider two special cases of $\sum_i |w_i| < \infty$ and $\sum_i |w_i| = \infty$. Generalizations are possible at the expense of additional details.

THEOREM 10. *Assume that* $\sum_{i=1}^{\infty} |w_i i^{1-\beta} L(i)| < \infty$, $E(\varepsilon_1^4) < \infty$ *and* $\sup_j \|K(X_0, X_j)\| < \infty$. *Then under the regularity conditions* (K1)–(K3), *we have* $n^{-1/2} \sum_{1 \leq i_1, i_2 \leq n} L_{i_1, i_2} \xrightarrow{d} \mathrm{N}(0, \sigma^2)$ *for some* $\sigma^2 < \infty$.

THEOREM 11. *Let* $w_i \equiv 1$, *and assume that* $E(\varepsilon_1^4) < \infty$, *and* $\sup_{i,j} \|K(X_i, X_j)\| < \infty$, $\sup_{i,j,\ell} \|K(\tilde{X}_i, \tilde{X}_j)\| < \infty$. *Then under the regularity conditions* (K1)–(K6), *we have* $n^{-3/2} \sum_{1 \leq i_1, i_2 \leq n} L_{i_1, i_2} \xrightarrow{d} \mathrm{N}(0, \sigma^2)$ *for some* $\sigma^2 < \infty$.

We conjecture that Theorems 10 and 11 can be made more general by dropping the restrictions on the $w_i$. While that generality is not achieved in this paper, the two theorems do already cover a wide range of interesting results. In particular, numerous limit theorems for the partial sum in the context of long-memory linear process [cf. Ho and Hsing (1997)] are special cases.

As explained earlier, the asymptotic distribution of $U_n$ is determined by one term or a combination of terms on the right-hand side of (28). The asymptotic behavior of $U_n$ is described by Theorems 10 and 11, while those for the "non-central" terms $Z_{n,r}$ are typically more straightforward but must be considered case by case. Let us take any $l_1, \ldots, l_r \in \{1, 2\}$ with $1 \leq r \leq \rho$ and $(r_1, r_2) = (p, q)$, and consider two special cases for the purpose of illustration.

First let us consider the case where $|w_t|$ is summable. Note that under general conditions we expect

$$K_{0,t,-\infty}(x_1, x_2) \to G(x_1, x_2) := EK(\hat{X}_1 + x_1, \hat{X}_2 + x_2),$$

where $\hat{X}_1, \hat{X}_2$ are i.i.d. that have the same distribution as $X_1$. Hence we assume that the $w_t D_{l_1,\ldots,l_r} K_{1,t,-\infty}(0,0)$ are absolutely summable in $t$ and

$$\lim_{n\to\infty}\left(w_0 D_{l_1,\ldots,l_r}K_{1,1,-\infty}(0,0) + 2\sum_{t=1}^{n} w_t D_{l_1,\ldots,l_r}K_{0,t,-\infty}(0,0)\right)$$
$$= C \in (-\infty, \infty).$$

Then it is not difficult to see [cf. Surgailis (1982) and Major (1980)] that

$$\left\|\sum_{1\leq i_1,i_2\leq n} w_{i_1-i_2} D_{l_1,\ldots,l_r} K_{i_1,i_2,-\infty}(0,0) \sum_{j_1>\cdots>j_r} \prod_{s=1}^{r} a_{i_{l_s}-j_s}\varepsilon_{j_s}\right\|^2$$

$$= \sum_{j_1>\cdots>j_r}\left(\sum_{1\leq i_1,i_2\leq n} w_{i_1-i_2} D_{l_1,\ldots,l_r} K_{i_1,i_2,-\infty}(0,0) \prod_{s=1}^{r} a_{i_{l_s}-j_s}\right)^2$$

$$\sim C^2 \sum_{j_1>\cdots>j_r}\left[\sum_{t=1}^{n}\prod_{s=1}^{r} a_{t-j_s}\right]^2$$

$$\sim C^2 n^{2-r(2\beta-1)} L^{2r}(n) \int_{u_1>\cdots>u_r}\left[\int_{x=0}^{1} \prod_{s=1}^{r}(x-u_s)_+^{-\beta}\,dx\right]^2 du_1\cdots du_r$$

and

$$n^{-1+r(\beta-1/2)}L^{-r}(n)$$

(37)
$$\times \sum_{1\leq i_1,i_2\leq n} w_{i_1-i_2} D_{l_1,\ldots,l_r} K_{i_1,i_2,-\infty}(0,0) \sum_{j_1>\cdots>j_r} \prod_{s=1}^{r} a_{i_{l_s}-j_s}\varepsilon_{j_s}$$

$$\xrightarrow{d} |C|\int_{u_1>\cdots>u_r}\left[\int_{x=0}^{1}\prod_{s=1}^{r}(x-u_s)_+^{-\beta}\,dx\right] dB(u_1)\cdots dB(u_r),$$

where the limit is expressed in the form of a multiple Wiener–Ito integral with $B$ denoting standard Brownian motion and $y_+ = \max(y, 0)$. Note that applying Theorem 10, the rate of $\sum_{i_1,i_2=1}^{n} L_{i_1,i_2}$ in this case is $n^{1/2}$, which is lower than that of $\sum_{1\leq i_1,i_2\leq n} D_{l_1,\ldots,l_r} H_{i_1,i_2,-\infty}(0,0) \sum_{j_1>\cdots>j_r} \prod_{s=1}^{r} a_{i_{l_s}-j_s}\varepsilon_{j_s}$.

As a second example, we consider an application in connection with Theorem 11 by assuming

(38) $\qquad w_i \equiv 1 \quad \text{and} \quad D_{l_1,\ldots,l_r}K_{0,t,-\infty}(0,0) \to C \in (-\infty, \infty).$

Under this assumption,

$$\left\|\sum_{1\leq i_1,i_2\leq n} D_{l_1,\ldots,l_r} K_{i_1,i_2,-\infty}(0,0) \sum_{j_1>\cdots>j_r} \prod_{s=1}^{r} a_{i_{l_s}-j_s}\varepsilon_{j_s}\right\|^2$$

$$\sim C^2 n^{4-r(2\beta-1)} L^{2r}(n)$$

$$\times \int_{u_1 > \cdots > u_r} \left[ \int_{x_1=0}^1 \int_{x_2=0}^1 \prod_{s=1}^r (x_{l_s} - u_s)_+^{-\beta} \, dx_1 \, dx_2 \right]^2 du_1 \cdots du_r$$

and

$$n^{-2+r(\beta-1/2)} L^{-r}(n) \sum_{1 \le i_1, i_2 \le n} D_{l_1,\ldots,l_r} K_{i_1,i_2,-\infty}(0,0) \sum_{j_1 > \cdots > j_r} \prod_{s=1}^r a_{i_{l_s}-j_s} \varepsilon_{j_s}$$

(39)
$$\xrightarrow{d} |C| \int_{u_1 > \cdots > u_r} \left[ \int_{x_1=0}^1 \int_{x_2=0}^1 (x_2 - x_1)^{-\beta} \right.$$

$$\left. \times \prod_{s=1}^r (x_{l_s} - u_s)_+^{-\beta} \, dx_1 \, dx_2 \right] dB(u_1) \cdots dB(u_r).$$

Now Theorem 11 implies that the rate of $\sum_{i_1,i_2=1}^n L_{i_1,i_2}$ is $n^{3/2}$, which is dominated by that of $\sum_{1 \le i_1,i_2 \le n} D_{i_1,\ldots,i_r} K_{i_1,i_2,-\infty}(0,0) \sum_{j_1 > \cdots > j_r} \prod_{s=1}^r a_{i_{l_s}-j_s} \varepsilon_{j_s}$.

In view of these examples and Theorems 10 and 11, in (28) one can refer to $\sum_{i_1,i_2=1}^n L_{i_1,i_2}$ the short-memory component, and $\sum_{i=1}^\rho Z_{n,r}$ the long-memory component of $U_n$.

Numerous applications result from these two simple cases. The following are some illustrations.

(a) Sample covariance function. Suppose

$$U_n = \sum_{i=1}^n T(X_i) T(X_{i+k}),$$

where $T$ is some function. $n^{-1} U_n$ is an estimator of $E[T(X_1) T(X_{1+k})]$. So $H_{i_1,i_2}(x_1,x_2) = w_{i_1-i_2} K(x_1,x_2)$, where $w_{|i_1-i_2|} = I(|i_1-i_2| = k)$ and $K(x_1,x_2) = T(x_1)T(x_2)$. Thus Theorem 10 applies, where the asymptotic distribution rests on $T$. Let us consider an example by assuming $k \ge 2$, $T(x) = x^2$ and $E(X_1 X_{1+k}^j) = 0, j = 1, 2$. It is easy to see that

$$w_t D_{l_1,\ldots,l_r} K_{0,t,-\infty}(0,0) = \begin{cases} 2E(X_1^2), & \text{if } t = k \text{ and} \\ & (l_1,\ldots,l_r) = (1,1) \text{ or } (2,2), \\ 0, & \text{otherwise.} \end{cases}$$

If $\beta \in (3/4, 1)$, then $\sum_{i=1}^\rho Z_{n,r} = 0$ and hence

$$\frac{1}{\sqrt{n}}(U_n - EU_n) \xrightarrow{d} N(0, \sigma^2).$$

If $\beta \in (1/2, 3/4)$, then $\sum_{i=1}^\rho Z_{n,r}$ is dominated by the term

$$Z_{n,2} = \sum_{1 \le i_1, i_2 \le n} (D_{1,1} + D_{2,2}) H_{i_1,i_2,-\infty}(0,0) \sum_{j_1 > \cdots > j_r} \prod_{s=1}^r a_{i_{l_s}-j_s} \varepsilon_{j_s}$$

$$= 8E(X_1^2) \sum_{i=1}^{n-k} \sum_{j_1 > j_r} a_{i-j_1} a_{i+k-j_2} \varepsilon_{j_1} \varepsilon_{j_2}.$$

Hence (37) and the discussion leading to it give

$$[8E(X_1^2)]^{-1} n^{-1+2(\beta-1/2)} L^{-2}(n)(U_n - EU_n)$$
$$\xrightarrow{d} \int_{u_1 > u_2} \left[ \int_{x=0}^1 [(x-u_1)_+(x-u_2)_+]^{-\beta} dx \right] dB(u_1) dB(u_2).$$

(b) $U$- and $V$-statistics. The asymptotic distribution of statistics of the form

$$\sum_{1 \le i_1 \ne i_2 \le n} K(X_{i_1}, X_{i_2}) \quad \text{or} \quad \sum_{1 \le i_1, i_2 \le n} K(X_{i_1}, X_{i_2})$$

can be considered using Theorems 10 and 11. Note that the partial-sum theory developed in Ho and Hsing (1997) is readily recovered here by letting $K(x_1, x_2) = (h(x_1) + h(x_2))/2$. We give another example here, the Wilcoxon one-sample statistic and the signed-rank statistic, for which the asymptotic distribution is not seen elsewhere.

Let

$$K_{i_1,i_2}(x_1, x_2) = I(x_1 + x_2 > 0).$$

Then $[n(n-1)]^{-1} \sum_{1 \le i_1 \ne i_2 \le n} K(X_{i_1}, X_{i_2})$ is called the Wilcoxon one-sample statistic. Let us assume for simplicity that $\varepsilon_j$ has a normal distribution and that the marginal of $\{X_i\}$ is standard normal. Simple calculations give

$$K_{i_1,i_2,-\infty}(x_1, x_2) = P(X_{i_1} + X_{i_2} > -(x_1 + x_2)) = 1 - \Phi\left(-\frac{x_1 + x_2}{\sqrt{2(1 + \rho_{i_1 - i_2})}}\right),$$

where $\rho_n = EX_1 X_{1+n}$ and $\phi$ is the standard normal p.d.f. Hence,

$$K_{i_1,i_2,-\infty}^{(1,0)}(0,0) = \frac{1}{\sqrt{2(1 + \rho_{i_1-i_2})}} \phi(0) \to \frac{1}{2\sqrt{\pi}} \quad \text{as } |i_1 - i_2| \to \infty.$$

With $C = 1/(2\sqrt{\pi})$ in (38), it follows from (39) that, for each $\beta \in (1/2, 1)$,

$$(40) \quad n^{\beta-5/2} L^{-1}(n)(U_n - EU_n) \xrightarrow{d} \frac{1}{\sqrt{\pi}} \int_{u \in \Re} \int_{x=0}^1 [(x-u)_+]^{-\beta} dx \, dB(u),$$

which has a normal distribution. A related statistics is the signed-rank statistic

$$W_n = \sum_{i=1}^n \Psi_i R_i^+,$$

where $\Psi_i$ = sign of $X_i$ and $R_i^+$ = the rank of $|X_i|$ among $|X_1|, \ldots, |X_n|$. It can be shown [cf. Randles and Wolfe (1979)] that $W_n = U_n(1 + o_p(1))$ and hence the asymptotic distribution can be derived in exactly the same way.

There are situations where $\sum_{i_1,i_2=1}^n L_{i_1,i_2}$ as well as $Z_{n,r}, 1 \le r \le \rho$, all equal zero. Then the asymptotic distribution will be determined by the lowest-order non-trivial $Z_{n,r}$. This is exemplified by certain $U$-statistics with degenerate kernels, kernels which satisfy $\int_{x_1} K(x_1, x_2) \, dF(x_1) = 0$ for all $x_2$. See Dehling and Taqqu (1989, 1991) and Ho and Hsing (1996). The approach in those references overlaps and complements the approach described here.

## 6. Proofs

PROOF OF THEOREM 1. Let
$$\xi_m(\mathbf{Z}_t) = \sum_{i=t-m}^{t+m} L_{t,i} \text{ and } S_n(\xi_m) = \sum_{t=1}^n \xi_m(\mathbf{Z}_t), \qquad m \ge 1.$$

Then (5) implies
$$\sum_{t=0}^\infty \|\mathcal{P}_0(\xi_m(\mathbf{Z}_t))\| \le \sum_{t=0}^\infty \sum_{i=t-m}^{t+m} |w_{t-i}|\theta_{t,i} \le \sum_{t=0}^\infty \sum_{k=-m}^m |w_k|\theta_{t,t-k} < \infty$$

which entails $S_n(\xi_m)/\sqrt{n} \Rightarrow N(0, \sigma_m^2)$ for some $\sigma_m^2 < \infty$ by Theorem 1 in Woodroofe (1992) [see also Lemma 5 in Wu (2003)].

Let LIM be $\limsup_{m \to \infty} \limsup_{n \to \infty}$. It remains to verify that

(41)
$$\text{LIM} \frac{1}{\sqrt{n}} \left\| \sum_{1 \le i_1, i_2 \le n} L_{i_1, i_2} - S_n(\xi_m) \right\| = \text{LIM} \frac{1}{\sqrt{n}} \left\| \sum_{|i_1 - i_2| > m, \, 1 \le i_1, i_2 \le n} L_{i_1, i_2} \right\| = 0.$$

To this end, note that the projections $\mathcal{P}_t$ are orthogonal and hence

(42)
$$\left\| \sum_{|i_1-i_2|>m,\, 1\le i_1,i_2 \le n} L_{i_1,i_2} \right\|^2 = \sum_{t=-\infty}^n \left\| \mathcal{P}_t \sum_{|i_1-i_2|>m,\, 1 \le i_1,i_2 \le n} L_{i_1,i_2} \right\|^2$$
$$\le \sum_{t=-\infty}^n \left[ \sum_{|i_1-i_2|>m,\, 1 \le i_1, i_2 \le n} \|\mathcal{P}_t L_{i_1,i_2}\| \right]^2$$
$$\le 2\left(\sum_{t=1}^n + \sum_{t=-\infty}^0\right)\left[\sum_{i_2=1}^n \sum_{i_1=i_2+m+1}^n \|\mathcal{P}_t L_{i_1,i_2}\|\right]^2.$$

Making use of the fact that $\mathcal{P}_t L_{i_1,i_2} = 0$ for $t \ge i_1 \vee i_2$,
$$\sum_{t=1}^n \left(\sum_{i_2=1}^n \sum_{i_1=i_2+m+1}^n \|\mathcal{P}_t L_{i_1,i_2}\|\right)^2 = \sum_{t=1}^n \left[\sum_{k=m+1}^{n-1} \sum_{i=t-k}^{n-k} \|\mathcal{P}_t L_{i+k,i}\|\right]^2$$

$$
\begin{aligned}
(43) \quad &= \sum_{t=1}^{n}\left[\sum_{k=m+1}^{n-1}\sum_{i=0}^{n-t}|w_k|\theta_{i,i-k}\right]^2 \\
&\leq n\left[\sum_{k=m+1}^{n-1}\sum_{i=0}^{n}|w_k|\theta_{i,i-k}\right]^2,
\end{aligned}
$$

and similarly,

$$
\sum_{t=-\infty}^{0}\left[\sum_{i_2=1}^{n}\sum_{i_1=i_2+m+1}^{n}\|\mathcal{P}_t L_{i_1,i_2}\|\right]^2
$$

$$
\begin{aligned}
(44) \quad &= \sum_{t=-\infty}^{0}\left[\sum_{k=m+1}^{n-1}\sum_{i=1}^{n-k}\|\mathcal{P}_t L_{i+k,i}\|\right]^2 \leq \sum_{t=0}^{\infty}\left[\sum_{k=m+1}^{n-1}\sum_{i=1}^{n-k}|w_k|\theta_{i+k+t,i+t}\right]^2 \\
&\leq C\sum_{t=0}^{\infty}\sum_{k=m+1}^{n-1}\sum_{i=1}^{n-k}|w_k|\theta_{i+k+t,i+t} \leq C\sum_{k=m+1}^{\infty}\sum_{i=1}^{n}\sum_{t=0}^{\infty}|w_k|\theta_{i+k+t,i+t} \\
&\leq C\sum_{k=m+1}^{\infty}\sum_{i=1}^{n}\sum_{j=0}^{\infty}|w_k|\theta_{j,j-k} = Cn\sum_{k=m+1}^{\infty}\sum_{j=0}^{\infty}|w_k|\theta_{j,j-k},
\end{aligned}
$$

in view of

$$
\sum_{k=m+1}^{n-1}\sum_{i=1}^{n-k}|w_k|\theta_{i+k+t,i+t} \leq \sum_{k=1}^{\infty}|w_k|\sum_{i=1}^{\infty}\theta_{i+k+t,i+t}
$$

$$
\leq \sum_{k=1}^{\infty}|w_k|\sum_{j=1}^{\infty}\theta_{j,j-k} =: C < \infty.
$$

Hence (41) follows from from (3), (42)–(44). □

PROOF OF THEOREM 2. By Cauchy's inequality, we have $\|\mathcal{P}_0(L_{i_1,i_2} - \tilde{L}_{i_1,i_2})\| \leq |w_{i_1-i_2}|\delta_\ell$. By the triangle and Cauchy's inequalities, we also have

$$
\|\mathcal{P}_0(L_{i_1,i_2} - \tilde{L}_{i_1,i_2})\| \leq C|w_{i_1-i_2}|\hat{\theta}_{i_1,i_2}.
$$

Thus there exists a $C > 0$ such that, for all $i_1, i_2, \ell$,

$$
(45) \quad \|\mathcal{P}_0(L_{i_1,i_2} - \tilde{L}_{i_1,i_2})\| \leq C|w_{i_1-i_2}|\min(\hat{\theta}_{i_1,i_2}, \delta_\ell).
$$

In the sequel let LIM stand for $\limsup_{\ell\to\infty}\limsup_{n\to\infty}$, and let $C$ stand for a constant which may vary from line to line. By the proof of Theorem 1, we have

$$
\Delta := \mathrm{LIM}\frac{1}{nW_n^2}\left\|\sum_{1\leq i_1,i_2\leq n}(L_{i_1,i_2} - \tilde{L}_{i_1,i_2})\right\|^2 \leq \mathrm{I} + \mathrm{II},
$$

where

$$\mathrm{I} = \mathrm{LIM} \frac{C}{nW_n^2} \sum_{t=1}^{n} \left[ \sum_{k=0}^{n-1} \sum_{i=0}^{n-t} |w_k| \min(\hat{\theta}_{i,i-k}, \delta_\ell) \right]^2,$$

$$\mathrm{II} = \mathrm{LIM} \frac{C}{nW_n^2} \sum_{t=0}^{\infty} \left[ \sum_{k=0}^{n-1} \sum_{i=1+k+t}^{n+t} |w_k| \min(\hat{\theta}_{i,i-k}, \delta_\ell) \right]^2.$$

By the assumptions, we have

$$\mathrm{I} \leq C \limsup_{\ell \to \infty} \sup_{k \geq 0} \sum_{i=0}^{\infty} \min(\hat{\theta}_{i,i-k}, \delta_\ell) = C \lim_{\delta \to 0} \sup_{k \geq 0} \sum_{i=0}^{\infty} \min(\hat{\theta}_{i,i-k}, \delta) = 0,$$

and

$$\mathrm{II} \leq \mathrm{LIM} \frac{C}{nW_n^2} \sum_{t=0}^{\infty} \left( \sum_{k=0}^{n-1} \sum_{i=1+t}^{n+t} |w_k| \min(\hat{\theta}_{i,i-k}, \delta_\ell) \right) \left( \sum_{k=0}^{n-1} \sum_{i=1}^{\infty} |w_k| \min(\hat{\theta}_{i,i-k}, \delta_\ell) \right)$$

$$\leq C \lim_{\delta \to 0} \left( \sup_{k \geq 0} \sum_{i=0}^{\infty} \min(\hat{\theta}_{i,i-k}, \delta) \right)^2 = 0.$$

Thus $\Delta = 0$ follows. $\square$

PROOF OF THEOREM 3. The plan of the proof is to show that, for every fixed $\ell \geq 1$,

(46) $$(nW_n^2)^{-1/2} \sum_{1 \leq i_1, i_2 \leq n} \tilde{L}_{i_1, i_2} \xrightarrow{d} \mathrm{N}(0, \tilde{\sigma}^2) \quad \text{as } n \to \infty$$

for some finite $\tilde{\sigma}^2$. It follows then from Theorem 2 that $\tilde{\sigma}^2$ is Cauchy in $\ell$ and hence converges to a finite constant as $\ell \to \infty$. By this and another application of Theorem 2, we conclude that (9) holds with $\sigma^2 = \lim_{\ell \to \infty} \tilde{\sigma}^2 \in (0, \infty)$. Thus, we will focus on proving (46) for a fixed $\ell$. Observe that

$$\left\| \sum_{1 \leq i_1, i_2 \leq n, |i_1 - i_2| \leq \ell} \tilde{L}_{i_1, i_2} \right\| = O\left( \sqrt{n} \sum_{k=0}^{\ell} |w_k| \right) = o[(nW_n^2)^{1/2}]$$

since $W_n \to \infty$. Thus it suffices to show that

(47) $$(nW_n^2)^{-1/2} \sum_{1 \leq i_1, i_2 \leq n, |i_1 - i_2| > \ell} \tilde{L}_{i_1, i_2} \xrightarrow{d} \mathrm{N}(0, \tilde{\sigma}^2) \quad \text{as } n \to \infty$$

for some finite $\tilde{\sigma}^2$. Fix $i_1, i_2$ with $|i_1 - i_2| > \ell$, and observe that $\tilde{\mathbf{Z}}_{i_1}$ and $\tilde{\mathbf{Z}}_{i_2}$ are i.i.d. Now define $\tilde{J}_{i_1, i_2}(\tilde{\mathbf{Z}}_{i_1}) = E[\tilde{L}_{i_1, i_2} | \tilde{\mathbf{Z}}_{i_1}]$ and

$$\tilde{R}_{i_1, i_2} = \tilde{L}_{i_1, i_2} - \tilde{J}_{i_1, i_2}(\tilde{\mathbf{Z}}_{i_1}) - \tilde{J}_{i_1, i_2}(\tilde{\mathbf{Z}}_{i_2}).$$

Let $i_1 = r_1 + q_1\ell$ and $i_2 = r_2 + q_2\ell$, where integers $0 \leq r_1, r_2 \leq \ell - 1$ and $0 \leq q_1 \neq q_2 \leq q = \lfloor n/\ell \rfloor$. Since $\tilde{R}_{r_1+q_1\ell, r_2+q_2\ell}$ are uncorrelated for different pairs $q_1 < q_2$,

$$\left\| \sum_{1 \leq q_1 < q_2 \leq q} \tilde{R}_{r_1+q_1\ell, r_2+q_2\ell} \right\|^2 \leq C \sum_{1 \leq q_1 < q_2 \leq q} w_{r_1+q_1\ell-r_2-q_2\ell}^2,$$

which in conjunction with Cauchy's inequality by summing over $r_1, r_2 = 0, \ldots, \ell - 1$ yields that

$$\left\| \sum_{1 \leq i_1, i_2 \leq n, |i_1-i_2| > \ell} \tilde{R}_{i_1,i_2} \right\|^2 \leq C \sum_{1 \leq r_1, r_2 \leq \ell} \sum_{1 \leq q_1 < q_2 \leq q} w_{r_1+q_1\ell-r_2-q_2\ell}^2$$

$$= C \sum_{1 \leq i_1, i_2 \leq n} w_{i_1-i_2}^2 \leq C \sum_{k=0}^{n} (n-k) w_k^2 = o(nW_n^2).$$

So (47) will follow from

(48) $\quad (nW_n^2)^{-1/2} \sum_{1 \leq i_1, i_2 \leq n, |i_1-i_2| > \ell} \tilde{J}_{i_1,i_2}(\tilde{\mathbf{Z}}_{i_1}) \xrightarrow{d} \mathrm{N}(0, \tilde{\sigma}^2) \quad$ as $n \to \infty$

for some finite $\tilde{\sigma}^2$, which can be easily shown to hold by the central limit theorem for $\ell$-dependent processes. $\square$

PROOF OF THEOREM 10. We verify (3). By Lemmas 12(iii) and 13, we have, for $k \geq 0$,

$$\sum_{i=k}^{\infty} \theta_{i,i-k} \leq C \sum_{i=k}^{\infty} \psi_{i-k} = C \sum_{i=0}^{\infty} \psi_i < \infty.$$

Next for $i \geq 0$ and $j \leq -1$, we have by (K3) that $\|\mathcal{P}_0 K(X_i, X_j)\| = O(|a_i|)$ uniformly. It is easily seen that, for the "well-structured" part $K(X_i, X_j) - Y_{i,j}$, we also have $\|\mathcal{P}_0[K(X_i, X_j) - Y_{i,j}]\| = O(|a_i|)$ uniformly. Hence, for $k \geq 0$,

$$\sum_{i=0}^{k-1} \theta_{i,i-k} \leq \sum_{i=0}^{k-1} |a_i| \leq Ck^{1-\beta} L(k)$$

by Karamata's theorem. Hence

$$\sum_{k=0}^{\infty} |w_k| \sum_{i=0}^{k-1} \theta_{i,i-k} \leq C \sum_{k=0}^{\infty} |w_k| k^{1-\beta} L(k) < \infty.$$

$\square$

PROOF OF THEOREM 11. For each $i_1 \geq i_2$, let

$$\mathcal{Y}_{i_1,i_2} = K(X_{i_1}, X_{i_2})$$
$$- \sum_{t=i_2+1}^{i_1} \sum_{r=1}^{\rho \wedge (t-i_2)} K_{i_1,i_2,i_2}^{(r,0)}(X_{i_1,-\infty,i_2}, X_{i_2}) \sum_{t=j_1>\cdots>j_r\geq i_2+1} \prod_{s=1}^{r} a_{i_1-j_s}\varepsilon_{j_s}$$
$$- \sum_{r=0}^{\rho} \sum_{|\mathbf{l}|=r} D_{\mathbf{l}} K_{i_1,i_2,-\infty}(0,0) \sum_{i_2 \geq j_1>\cdots>j_r} \prod_{s=1}^{r} a_{i_{l_s}-j_s}\varepsilon_{j_s}$$

and $\mathcal{Y}_{i_1,i_2} = \mathcal{Y}_{i_2,i_1}$ if $i_1 < i_2$. We first apply Theorem 3 to show that

$$(49) \qquad \frac{1}{n^{3/2}} \sum_{i_1,i_2=1}^{n} \mathcal{Y}_{i_1,i_2} \xrightarrow{d} \mathrm{N}(0, \sigma^2)$$

for some $\sigma^2 < \infty$. For the process $\{\mathcal{Y}_{i,j}\}$, it follows easily from (36) that $\lim_{\ell \to \infty} \delta_\ell = 0$. Hence we focus on verifying the condition (7). Observe that for $i_1 \geq i_2$, $\mathcal{P}_t \mathcal{Y}_{i_1,i_2}$ is equal to

$$\mathcal{P}_t K(X_{i_1}, X_{i_2}) - \sum_{r=1}^{\rho \wedge (t-i_2)} H_{i_1,i_2,i_2}^{(r,0)}(X_{i_1,-\infty,i_2}, X_{i_2}) \sum_{t=j_1>\cdots>j_r\geq i_2+1} \prod_{s=1}^{r} a_{i_1-j_s}\varepsilon_{j_s}$$

if $i_2 < t \leq i_1$, and is equal to $\mathcal{P}_t Y_{i_1,i_2}$ if $t \leq i_2$. Hence by Lemmas 12 and 14, for $i_1 \geq \max(0, i_2)$,

$$\hat{\theta}_{i_1,i_2} = \sup_{\ell \geq 0} \mathcal{P}_0 \tilde{\mathcal{Y}}_{i_1,i_2} \leq \begin{cases} C\xi_{i_1,i_2}, & i_2 < 0 \leq i_1, \\ C\psi_{i_2}, & i_2 \geq 0. \end{cases}$$

Note that $\xi_{i_1,i_2} = \psi_{i_1}$ if $i_2 \leq -\rho$. Thus

$$\sum_{i=0}^{\infty} \min(\hat{\theta}_{i,i-k}, \epsilon) = \sum_{i=0}^{k-1} \min(\hat{\theta}_{i,i-k}, \epsilon) + \sum_{i=0}^{\infty} \min(\hat{\theta}_{i+k,i}, \epsilon)$$
$$\leq C\left(\epsilon + \sum_{i=0}^{\infty} \min(\psi_i, \epsilon)\right) \to 0 \qquad \text{as } \epsilon \to 0,$$

since the $\psi_i$ are summable by Lemma 13. Hence the condition (7) is proved for $\{\mathcal{Y}_{i,j}\}$ and the proof for (49) is complete.

Next observe that

$$Y_{i_1,i_2} = \mathcal{Y}_{i_1,i_2} + W_{i_1,i_2}$$
$$:= \mathcal{Y}_{i_1,i_2} + \sum_{r=1}^{\rho} R_{i_1,i_2,r} \sum_{t=i_2+r}^{i_1} \sum_{t=j_1>\cdots>j_r\geq i_2+1} \prod_{s=1}^{r} a_{i_1-j_s}\varepsilon_{j_s},$$

where
$$R_{i_1,i_2,r} = K^{(r,0)}_{i_1,i_2,i_2}(X_{i_1,-\infty,i_2}, X_{i_2})$$
$$- \sum_{r'=0}^{\rho-r} \sum_{|\mathbf{l}|=r'} D_{\mathbf{l}} K^{(r,0)}_{i_1,i_2,-\infty}(0,0) \sum_{i_2 \geq j_1 > \cdots > j_{r'}} \prod_{s=1}^{r'} a_{i_{l_s}-j_s} \varepsilon_{j_s}.$$

Hence the conclusion of the theorem follows from (49) and Lemma 15. □

For $\mathbf{k} \in \mathcal{K}$ with $|\mathbf{k}| \leq \rho$, define
$$M^{(\mathbf{k})}_{i_1,i_2} = D_{\mathbf{k}} K_{i_1,i_2,0}(X_{i_1,-\infty,0}, X_{i_2,-\infty,0})$$
$$- \sum_{r=0}^{\rho-|\mathbf{k}|} \sum_{|\mathbf{l}|=r} D_{\mathbf{l}} D_{\mathbf{k}} K_{i_1,i_2,-\infty}(0,0) \sum_{0 \geq j_1 > j_2 > \cdots > j_r} \prod_{s=1}^{r} a_{i_{l_s}-j_s} \varepsilon_{j_s}.$$

LEMMA 12. *Assume that $E(\varepsilon_1^4) < \infty$, and (K1) and (K2) hold. Then there exists a constant $C$, independent of $i_1, i_2 \geq n_0$, such that*

(i) *for all $\mathbf{k} \in \mathcal{K}$ with $|\mathbf{k}| \leq \rho - 1$,*
$$\|\mathcal{P}_0 M^{(\mathbf{k})}_{i_1,i_2}\|^2 \leq C[a_{i_1}^4 + a_{i_2}^4 + a_{i_1}^2 \|M^{(1,\mathbf{k})}_{i_1,i_2}\|^2 + a_{i_2}^2 \|M^{(2,\mathbf{k})}_{i_1+1,i_2+1}\|^2];$$

(ii) *for all $\mathbf{k} \in \mathcal{K}$ such that $|\mathbf{k}| \leq \rho$,*
$$\|M^{(\mathbf{k})}_{i_1,i_2}\|^2 \leq C[A_{i_1}(4) + A_{i_2}(4) + A^{\rho-|\mathbf{k}|+1}_{i_1}(2) + A^{\rho-|\mathbf{k}|+1}_{i_2}(2)];$$

(iii) $\|\mathcal{P}_0 Y_{i_1,i_2}\|^2 \leq C(\psi_{i_1}^2 + \psi_{i_2}^2)$, *where $\psi_i = |a_i|[|a_i| + \sqrt{A_{i+1}(4) + A^\rho_{i+1}(2)}]$.*

PROOF. We first prove (i). Define
$$B_{\mathbf{k}}(i_1,i_2) = \sum_{r=1}^{\rho-|\mathbf{k}|} \sum_{|\mathbf{l}|=r} D_{\mathbf{l}} D_{\mathbf{k}} K_{i_1,i_2,-\infty}(0,0) \sum_{0=j_1>j_2>\cdots>j_r} \prod_{s=1}^{r} a_{i_{l_s}-j_s}\varepsilon_{j_s},$$
$$B_{1,\mathbf{k}}(i_1,i_2) = \sum_{r=1}^{\rho-|\mathbf{k}|} \sum_{|\mathbf{l}|=r-1} D_{\mathbf{l}} D_{1,\mathbf{k}} K_{i_1,i_2,-\infty}(0,0) \sum_{-1\geq j_2>\cdots>j_r} \prod_{s=2}^{r} a_{i_{l_s}-j_s}\varepsilon_{j_s},$$
$$B_{2,\mathbf{k}}(i_1,i_2) = \sum_{r=1}^{\rho-|\mathbf{k}|} \sum_{|\mathbf{l}|=r-1} D_{\mathbf{l}} D_{2,\mathbf{k}} K_{i_1,i_2,-\infty}(0,0) \sum_{-1\geq j_2>\cdots>j_r} \prod_{s=2}^{r} a_{i_{l_s}-j_s}\varepsilon_{j_s}.$$

Then $B_{\mathbf{k}}(i_1,i_2) = a_{i_1}\varepsilon_0 B_{1,\mathbf{k}}(i_1,i_2) + a_{i_2}\varepsilon_0 B_{2,\mathbf{k}}(i_1,i_2)$. Observe that, for $\iota = 1$ and 2,
$$\|D_{\iota,\mathbf{k}} K_{i_1,i_2,-1}(X_{i_1,-\infty,-1}, X_{i_2,-\infty,-1}) - B_{\iota,\mathbf{k}}(i_1,i_2)\| = \|M^{(\iota,\mathbf{k})}_{i_1+1,i_2+1}\|.$$

Hence by the triangle inequality,

$$\|\langle \nabla D_{\mathbf{k}} K_{i_1,i_2,-1}(X_{i_1,-\infty,-1}, X_{i_2,-\infty,-1}), (a_{i_1}\varepsilon_0, a_{i_2}\varepsilon_0)\rangle - B_{\mathbf{k}}(i_1, i_2)\|$$

$$(50) \quad \leq \sum_{\iota=1}^{2} |a_{i_\iota}| \|D_{\iota,\mathbf{k}} K_{i_1,i_2,-1}(X_{i_1,-\infty,-1}, X_{i_2,-\infty,-1}) - B_{\iota,\mathbf{k}}(i_1, i_2)\|$$

$$= \sum_{\iota=1}^{2} |a_{i_\iota}| \|M_{i_1+1,i_2+1}^{(\iota,\mathbf{k})}\|.$$

By (29) and the triangle inequality,

$$\|\mathcal{P}_0 M_{i_1,i_2}^{(\mathbf{k})}\|$$
$$\leq \|D_{\mathbf{k}} K_{i_1,i_2,0}(X_{i_1,-\infty,0}, X_{i_2,-\infty,0})$$
$$\quad - D_{\mathbf{k}} K_{i_1,i_2,-1}(X_{i_1,-\infty,-1}, X_{i_2,-\infty,-1})$$
$$\quad - \langle \nabla D_{\mathbf{k}} K_{i_1,i_2,-1}(X_{i_1,-\infty,-1}, X_{i_2,-\infty,-1}), (a_{i_1}\varepsilon_0, a_{i_2}\varepsilon_0)\rangle\|$$
$$\quad + \|\langle \nabla D_{\mathbf{k}} K_{i_1,i_2,-1}(X_{i_1,-\infty,-1}, X_{i_2,-\infty,-1}), (a_{i_1}\varepsilon_0, a_{i_2}\varepsilon_0)\rangle - B_{\mathbf{k}}(i_1, i_2)\|,$$

from which (i) follows in view of (30) and (50).

To establish (ii) we will adopt a backward induction argument. First, for $|\mathbf{k}| = \rho$, since $M_{i_1,i_2}^{(\mathbf{k})} = D_{\mathbf{k}} K_{i_1,i_2,0}(X_{i_1,-\infty,0}, X_{i_2,-\infty,0}) - D_{\mathbf{k}} K_{i_1,i_2,-\infty}(0,0)$, (ii) follows from (31). Next we make the induction assumption that (ii) holds for all $\mathbf{k}$ with $|\mathbf{k}| = m \geq 1$ and we wish to show that it holds for all $|\mathbf{k}| = m - 1$. By (i) we have, for any $\mathbf{k}$ with $|\mathbf{k}| = m - 1$,

$$\|\mathcal{P}_0 M_{i_1,i_2}^{(\mathbf{k})}\|^2 \leq C[a_{i_1}^4 + a_{i_2}^4 + a_{i_1}^2 \|M_{i_1,i_2}^{(1,\mathbf{k})}\|^2 + a_{i_2}^2 \|M_{i_1+1,i_2+1}^{(2,\mathbf{k})}\|^2].$$

Since the projections $\mathcal{P}_t$ are orthogonal and $\mathcal{P}_t M_{i_1,i_2}^{(\mathbf{k})} = 0$ for $t \geq 1$,

$$\|M_{i_1,i_2}^{(\mathbf{k})}\|^2 = \sum_{t=-\infty}^{0} \|\mathcal{P}_t M_{i_1,i_2}^{(\mathbf{k})}\|^2 = \sum_{t=0}^{\infty} \|\mathcal{P}_0 M_{i_1+t,i_2+t}^{(\mathbf{k})}\|^2$$

$$\leq C\bigg([A_{i_1}(4) + A_{i_2}(4)] + \sum_{t=0}^{\infty}(a_{i_1+t}^2 + a_{i_2+t}^2)[A_{i_1+t}(4) + A_{i_2+t}(4)]$$

$$+ \sum_{t=0}^{\infty}(a_{i_1+t}^2 + a_{i_2+t}^2)[A_{i_1+t}^{\rho-m+1}(2) + A_{i_2+t}^{\rho-m+1}(2)]\bigg)$$

by the induction assumption. Now the induction is complete since $\sum_{t=0}^{\infty} a_{i+t}^2 \times A_{i+t}(4) = o[A_i(4)]$ and $\sum_{t=0}^{\infty} a_{i+t}^2 A_{i+t}^{\rho-m+1}(2) = O[A_i^{\rho-m+2}(2)]$.

Finally, (iii) follows readily from parts (i) and (ii) of this lemma by noting that $\mathcal{P}_t Y_{i_1,i_2}$ and $\mathcal{P}_0 Y_{i_1-t,i_2-t} = \mathcal{P}_0 M_{i_1-t,i_2-t}^{(0)}$ have the same distribution. $\square$

LEMMA 13. *If $\rho$ and $L(\cdot)$ satisfy (27), then $\sum_{i=1}^{\infty} \psi_i < \infty$.*

PROOF. By Karamata's theorem,

(51) $$A_n(k) = O(na_n^k) \quad \text{for } k \geq 2.$$

So the lemma easily follows from (27) after elementary calculations. □

For all $i_1 \geq 0 \geq i_2$, define

$$N_{i_1,i_2}^{(k)} = \sum_{t=i_2+1}^{0} \mathcal{P}_t N_{i_1,i_2}^{(k)},$$

where $\mathcal{P}_t N_{i_1,i_2}^{(k)}$ is equal to

$$\mathcal{P}_t K_{i_1,i_2,t}^{(k,0)}(X_{i_1,-\infty,t}, X_{i_2})$$
$$- \sum_{r=1}^{(\rho \wedge (t-i_2))-k} K_{i_1,i_2,i_2}^{(k+r,0)}(X_{i_1,-\infty,i_2}, X_{i_2}) \sum_{t=j_1>\cdots>j_r \geq i_2+1} \prod_{s=1}^{r} a_{i_1-j_s} \varepsilon_{j_s}$$

if $0 \leq k \leq \rho \wedge (t - i_2)$ and 0 otherwise. The following lemma is very similar to Lemma 12.

LEMMA 14. *Assume that $E(\varepsilon_1^4) < \infty$ and conditions* (K4) *and* (K5) *hold. Let $i_1 \geq n_0 > 0 > i_2$ and write $\rho' = \rho \wedge (-i_2)$. Then there exists a constant $C$, independent of $i_1, i_2$, such that*

(i) *for all $0 \leq k \leq \rho' - 1$,*

$$\|\mathcal{P}_0 N_{i_1,i_2}^{(k)}\|^2 \leq C[a_{i_1}^4 + a_{i_1}^2 \|N_{i_1+1,i_2+1}^{(k+1)}\|^2];$$

(ii) *for all $0 \leq k \leq \rho'$,*

$$\|N_{i_1,i_2}^{(k)}\|^2 \leq C[A_{i_1}(4) + A_{i_1}^{\rho'-|\mathbf{k}|+1}(2)];$$

(iii) $\|\mathcal{P}_0 \mathcal{Y}_{i_1,i_2}\|^2 \leq C\xi_{i_1,i_2}^2$, *where* $\xi_{i_1,i_2} = |a_{i_1}|[|a_{i_1}| + \sqrt{A_{i_1+1}(4) + A_{i_1+1}^{\rho'}(2)}]$.

PROOF. We first prove (i). Fix $k \leq \rho' - 1$. By the triangle inequality and (33),

$$\|\mathcal{P}_0 N_{i_1,i_2}^{(k)}\| \leq \|K_{i_1,i_2,0}^{(k,0)}(X_{i_1,-\infty,0}, X_{i_2}) - K_{i_1,i_2,-1}^{(k,0)}(X_{i_1,-\infty,-1}, X_{i_2})$$
$$- a_{i_1} \varepsilon_0 K_{i_1,i_2,-1}^{(k+1,0)}(X_{i_1,-\infty,-1}, X_{i_2})\|$$
$$+ \left\| a_{i_1} \varepsilon_0 K_{i_1,i_2,-1}^{(k+1,0)}(X_{i_1,-\infty,-1}, X_{i_2}) \right.$$
$$\left. - \sum_{r=1}^{\rho'-k} K_{i_1,i_2,i_2}^{(k+r,0)}(X_{i_1,-\infty,i_2}, X_{i_2}) \sum_{0=j_1>j_2>\cdots>j_r \geq i_2+1} \prod_{s=1}^{r} a_{i_1-j_s} \varepsilon_{j_s} \right\|.$$

The first term on the right-hand side is bounded by $C|w_{i_1-i_2}|a_{i_1}^2$ by the assumption (34). The second term on the right-hand side is bounded by

$$|a_{i_1}| \bigg\| K^{(k+1,0)}_{i_1,i_2,-1}(X_{i_1,-\infty,-1}, X_{i_2}) - K^{(k+1,0)}_{i_1,i_2,i_2}(X_{i_1,-\infty,i_2}, X_{i_2})$$
$$- \sum_{r=1}^{\rho'-k-1} K^{(k+1+r,0)}_{i_1,i_2,i_2}(X_{i_1,-\infty,i_2}, X_{i_2}) \sum_{-1 \geq j_1 > \cdots > j_r \geq i_2+1} \prod_{s=1}^{r} a_{i_1-j_s}\varepsilon_{j_s} \bigg\|$$
$$\leq |a_{i_1}| \|N^{(k+1)}_{i_1+1,i_2+1}\| + |a_{i_1}| \|W_{i_1,i_2,k} - N^{(k+1)}_{i_1+1,i_2+1}\|,$$

where

$$W_{i_1,i_2,k} = K^{(k+1,0)}_{i_1+1,i_2+1,0}(X_{i_1+1,-\infty,0}, X_{i_2+1})$$
$$- K^{(k+1,0)}_{i_1+1,i_2+1,i_2+1}(X_{i_1+1,-\infty,i_2+1}, X_{i_2+1})$$
$$- \sum_{r=1}^{\rho'-k-1} K^{(k+1+r,0)}_{i_1+1,i_2+1,i_2+1}(X_{i_1+1,-\infty,i_2+1}, X_{i_2+1})$$
$$\times \sum_{0 \geq j_1 > \cdots > j_r \geq i_2+2} \prod_{s=1}^{r} a_{i_1+1-j_s}\varepsilon_{j_s}.$$

Observe that

$$\|W_{i_1,i_2,k} - N^{(k+1)}_{i_1+1,i_2+1}\|$$
$$\leq \sum_{t=i_2+1}^{0} \sum_{r=(\rho \wedge (t-i_2))-k+1}^{\rho'-k-1} \bigg\| K^{(k+1+r,0)}_{i_1+1,i_2+1,i_2+1}(X_{i_1+1,-\infty,i_2+1}, X_{i_2+1})$$
$$\times \sum_{t=j_1 > \cdots > j_r \geq i_2+2} \prod_{s=1}^{r} a_{i_1+1-j_s}\varepsilon_{j_s} \bigg\|$$
$$\leq C|a_{i_1}|.$$

Hence (i) follows.

To establish (ii), we will adopt a backward induction argument. First, for $k = \rho'$, (ii) follows from (35) in view of

$$N^{(\rho')}_{i_1,i_2} = K^{(\rho',0)}_{i_1,i_2,0}(X_{i_1,-\infty,0}, X_{i_2}) - K^{(\rho',0)}_{i_1,i_2,i_2}(X_{i_1,-\infty,i_2}, X_{i_2}).$$

Next we make the induction assumption that (ii) holds for $k = m \geq 1$ and we wish to show that it holds for $k = m-1$. By (i) we have

$$\|\mathcal{P}_0 N^{(m-1)}_{i_1,i_2}\|^2 \leq C[a_{i_1}^4 + a_{i_1}^2 \|N^{(m)}_{i_1,i_2}\|^2].$$

Since the projections $\mathcal{P}_t$ are orthogonal,

$$\|N_{i_1,i_2}^{(m)}\|^2 = \sum_{t=i_2+1}^{0} \|\mathcal{P}_t N_{i_1,i_2}^{(m)}\|^2 = \sum_{t=0}^{-i_2-1} \|\mathcal{P}_0 N_{i_1+t,i_2+t}^{(m)}\|^2$$

$$\leq Cw_{|i-j|}^2\left(A_{i_1}(4) + \sum_{t=0}^{\infty} a_{i_1+t}^2 A_{i_1+t}(4) + \sum_{t=0}^{\infty} a_{i_1+t}^2 A_{i_1+t}^{\rho'-m+1}(2)\right)$$

by the induction assumption. Now the induction is complete since $\sum_{t=0}^{\infty} a_{i+t}^2 \times A_{i+t}(4) = o[A_i(4)]$ and $\sum_{t=0}^{\infty} a_{i+t}^2 A_{i+t}^{\rho'-m+1}(2) = O[A_i^{\rho'-m+2}(2)]$.

Finally, (iii) follows readily from parts (i) and (ii) of this lemma by noting that $\mathcal{P}_t \mathcal{Y}_{i_1,i_2}$ and $\mathcal{P}_0 \mathcal{Y}_{i_1-t,i_2-t} = \mathcal{P}_0 N_{i_1-t,i_2-t}^{(0)}$ have the same distribution. $\square$

LEMMA 15. *Under the conditions of Theorem* 11,
$$n^{-3/2} \sum_{i,j=1}^{n} W_{i_1,i_2} \xrightarrow{p} 0.$$

PROOF. Recall that
$$R_{i_1,i_2,r} = K_{i_1,i_2,i_2}^{(r,0)}(X_{i_1,-\infty,i_2}, X_{i_2})$$
$$- \sum_{r'=0}^{\rho-r} \sum_{|\mathbf{l}|=r'} D_{\mathbf{l}} K_{i_1,i_2,-\infty}^{(r,0)}(0,0) \sum_{i_2 \geq j_1 > \cdots > j_{r'}} \prod_{s=1}^{r'} a_{i_{l_s}-j_s} \varepsilon_{j_s},$$

and hence we have for $i_1 \geq i_2 \geq 0$,
$$\|\mathcal{P}_0 R_{i_1,i_2,r}\| = \|\mathcal{P}_0 \mathcal{M}_{i_1,i_2}^{(r,0)}\| \leq C|a_{i_2}|[|a_{i_2}| + \sqrt{A_{i_2+1}(4) + A_{i_2+1}^{\rho-r}(2)}]$$
$$\sim Ca_{i_2}(ia_{i_2})^{(\rho-r)/2} \sim i_2^{-\delta} L^{\rho-r+1}(i_2) =: \eta_{i_2},$$

where $\delta = (\rho - r + 1)(\beta - 1/2) + 1/2$ by Lemma 12 and (51). Projecting iteratively, we obtain

$$\left\| \sum_{i_1 \geq i_2 = 1}^{n} R_{i_1,i_2,r} \sum_{i_2+r \leq t \leq i_1} \sum_{t=j_1 > \cdots > j_r \geq i_2+1} a_{i_1-j_1} \prod_{s=1}^{r} a_{i_1-j_s} \varepsilon_{j_s} \right\|^2$$

$$= \sum_{n \geq t_1 > \cdots > t_r \geq 1} \left\| \sum_{i_1=t_1}^{n} \sum_{i_2=1}^{t_r-1} \prod_{s=1}^{r} a_{i_1-t_s} R_{i_1,i_2,r} \right\|^2$$

$$= \sum_{n \geq t_1 > \cdots > t_r \geq 1} \sum_{t'=-\infty}^{t_r-1} \left\| \sum_{i_1=t_1}^{n} \sum_{i_2=1 \vee t'}^{t_r-1} \prod_{s=1}^{r} a_{i_1-t_s} \mathcal{P}_{t'} R_{i_1,i_2,r} \right\|^2$$

$$\leq C \sum_{n \geq t_1 > \cdots > t_r \geq 1} \sum_{t'=-\infty}^{t_r-1} \left( \sum_{i_1=t_1}^{n} \sum_{i_2=1 \vee t'}^{t_r-1} \prod_{s=1}^{r} a_{i_1-t_s} \eta_{i_2-t'} \right)^2.$$

Now approximating summations by integrals, the last expression can be seen to be asymptotically equal to $Cn^{r+5-2(r\beta+\delta)}(L(n))^{2(\rho-r+1)} \times INT$ for large $n$, where

$$INT = \int_{1>t_1>\cdots>t_r>0} \int_{t'=-\infty}^{t_r} \left(\int_{x=t_1}^1 \prod_{s=1}^r (x-t_s)^{-\beta}\right)^2 \left(\int_{y=0\vee t'}^{t_r-1} (y-t')^{-\delta}\right)^2.$$

Since

$$r+5-2(r\beta+\delta) = 4-(\rho+1)(2\beta-1) \le 3,$$

where the equality holds only if $L(n) \to 0$ by (27), the result will follow if we can show that $INT < \infty$, which is what we will do. First,

$$\int_{x=t_1}^1 \prod_{s=1}^r (x-t_s)^{-\beta} dx$$

$$\le (t_1-t_2)^{-(r\beta-1)} \int_0^\infty x^{-\beta}(x+1)^{-\beta} \prod_{s=3}^r \left(x + \frac{t_1-t_s}{t_1-t_2}\right)^{-\beta} dx$$

$$\le (t_1-t_2)^{-(r\beta-1)} \prod_{s=3}^r \left(\frac{t_1-t_s}{t_1-t_2}\right)^{-\beta} \int_0^\infty x^{-\beta}(x+1)^{-\beta} dx$$

$$\le C(t_1-t_2)^{-(2\beta-1)} \prod_{s=3}^r (t_1-t_s)^{-\beta}.$$

Hence we have

$$INT \le C \int_{1>t_1>\cdots>t_r>0} (t_1-t_2)^{-2(2\beta-1)} \prod_{s=3}^r (t_1-t_s)^{-2\beta}$$

$$\times \int_{t'=-\infty}^{t_r} \left(\int_{y=0\vee t'}^{t_r-1} (y-t')^{-\delta}\right)^2.$$

Writing

$$\int_{t'=-\infty}^{t_r} \left(\int_{y=0\vee t'}^{t_r-1} (y-t')^{-\delta} dy\right)^2 dt'$$

$$= \left(\int_{t'=-\infty}^{-1} + \int_{t'=-1}^0 + \int_{t'=0}^{t_r}\right) \left(\int_{y=0\vee t'}^{t_r-1} (y-t')^{-\delta} dy\right)^2 dt',$$

it is easy to see that all three integrals are uniformly bounded since $1/2 < \delta < 1$. Also, integrating iteratively from $t_r$ to $t_3$, we obtain

$$\int_{t_2>\cdots>t_r>0} \prod_{s=3}^r (t_1-t_s)^{-2\beta} dt_3\cdots dt_r \le C(t_1-t_2)^{-(r-2)(2\beta-1)}.$$

Since $r(2\beta - 1) < 1$,

$$\int_{1>t_1>\cdots>t_r>0}(t_1 - t_2)^{-2(2\beta-1)}\prod_{s=3}^{r}(t_1 - t_s)^{-2\beta}dt_1\cdots dt_r$$
$$\leq C\int_{t_1=0}^{1}\int_{t_2=0}^{t_1}(t_1 - t_2)^{-r(2\beta-1)}\,dt_1\,dt_2 < \infty.$$

Hence we conclude that $INT < \infty$ and the proof is complete. $\square$

**Acknowledgment.** The authors gratefully acknowledge the helpful and constructive comments provided by the anonymous referee.

Department of Statistics  
Texas A&M University  
College Station, Texas 77843  
USA  
e-mail: thsing@stat.tamu.edu

Department of Statistics  
University of Chicago  
5734 S. University Avenue  
Chicago, Illinois 60637  
USA  
e-mail: wbwu@galton.uchicago.edu